\documentclass[11pt,a4paper,reqno]{amsart}
\usepackage{amsaddr}
\usepackage{amsmath,amssymb,amsfonts,epsfig,mathrsfs,cite}
\usepackage[T1]{fontenc}
\usepackage{color}
\usepackage{array}
\usepackage{amsthm}
\usepackage{amstext}
\usepackage{graphicx}
\usepackage{setspace}
\usepackage{tcolorbox}
\usepackage{float}
\usepackage{subfigure}
\usepackage{hyperref} 
\usepackage[title]{appendix}
\usepackage{booktabs}
\usepackage{multirow}
\usepackage{algorithm}
\usepackage{algpseudocode}

\makeatletter
\@namedef{subjclassname@2020}{%
  \textup{2020} Mathematics Subject Classification}
\makeatother

\usepackage[margin=2.5cm]{geometry}
\usepackage{color}
\usepackage{enumitem}
\usepackage{amscd,psfrag}
\usepackage{yhmath}
\usepackage[mathscr]{eucal}

\usepackage{comment}

\allowdisplaybreaks[4]

\usepackage{slashed}

\makeatletter
\pdfpageheight\paperheight
\pdfpagewidth\paperwidth

\usepackage{epstopdf}
\usepackage{chngcntr}
\counterwithin{figure}{section}
\counterwithin{table}{section}
\counterwithout{equation}{section}
\usepackage{mathrsfs}

\usepackage{indentfirst}	

\usepackage[normalem]{ulem}
\theoremstyle{plain}

\newtheorem{definition}{Definition}[section]
\newtheorem{theorem}[definition]{Theorem}
\newtheorem*{theorem*}{Theorem}

\newtheorem{remark}[definition]{Remark}

\newtheorem*{remark*}{Remark}
\newtheorem*{sideremark*}{Side Remark}

\newtheorem*{claim*}{Claim}
\newtheorem*{q*}{Question}
\newtheorem{lemma}[definition]{Lemma}
\newtheorem{corollary}[definition]{Corollary}
\newtheorem*{corollary*}{Corollary}
\newtheorem{example}[definition]{Example}
\newtheorem{question}[definition]{Question}
\newtheorem{proposition}[definition]{Proposition}

\newcommand{\R}{\mathbb{R}}

\newcommand{\emb}{\hookrightarrow}

\newcommand{\p}{\partial}

\newcommand{\e}{\varepsilon}
\newcommand{\C}{\mathbb{C}}
\newcommand{\dd}{{\rm d}}
\newcommand{\g}{{\mathfrak{g}}}

\newcommand{\M}{{\mathcal{M}}}
\newcommand{\bra}{\left\langle}
\newcommand{\ket}{\right\rangle}
\newcommand{\bbra}{\bra\bra}
\newcommand{\kket}{\ket\ket}

\newcommand{\2}{{\sqrt{-1}}}

\newcommand{\E}{{\mathscr{E}}}
\newcommand{\hatE}{\widehat{\E}}

\newcommand{\lin}{{\mathcal{L}}}

\newcommand{\su}{{\mathfrak{su}}}
\newcommand{\so}{{\mathfrak{so}}}
\newcommand{\mt}{{M_{(\theta)}}}

\newcommand{\word}{{\bf Word}}
\newcommand{\wordn}{\word(n)}

\newcommand{\mw}{{\M[W]}}
\newcommand{\mwt}{{\mw_{(\theta)}}}

\newcommand{\tens}{{\T\left(\R^d\right)}}
\newcommand{\T}{{\mathcal{T}}}

\newcommand{\dset}{{\{1,\ldots,d\}}}
\newcommand{\nset}{{\{1,\ldots,n\}}}
\newcommand{\kset}{{\{1,\ldots,k\}}}
\newcommand{\tilm}{{\widetilde{M}}}
\newcommand{\matn}{{\mathfrak{gl}(n,\C)}}
\newcommand{\matk}{{\mathfrak{gl}(k,\C)}}
\newcommand{\cj}{{\mathcal{C}_J}}
\newcommand{\cw}{{\mathcal{C}_W}}
\newcommand{\ww}{{\overline{w}}}
\newcommand{\WW}{{\overline{W}}}

\newcommand{\dev}{{\mathscr{D}}}

\def\XXint#1#2#3{{\setbox0=\hbox{$#1{#2#3}{\int}$ }
\vcenter{\hbox{$#2#3$ }}\kern-.6\wd0}}

\newcommand{\art}{{\mathcal{A}_{\mathfrak{t},R}}}

\newcommand{\clart}{{\bf Clos}\left[\art\right]}

\title{Restricted Path Characteristic Function Determines the Law of Stochastic Processes}
\author{Siran Li$^a$, Zijiu Lyu$^b$, Hao Ni$^{c, \ast}$, Jiajie Tao$^c$}
\address{\small
$^a$ School of Mathematical Sciences $\&$ CMA-Shanghai, Shanghai Jiao Tong University, No.~6 Natural Science Buildings, 800 Dongchuan Road, Minhang District, Shanghai 200240, China \\
$^b$ Mathematical Institute, University of Oxford, Andrew Wiles Building, Woodstock Rd, Oxford OX2 6GG, UK \\
$^c$ Department of Mathematics, University College London, 25 Gordon St, London WC1H 0AY, UK \\
\medskip
\raggedright
Email addresses and ORCID iDs: 
\emph{
\begin{itemize}
    \item SIRAN LI: \href{mailto:siran.li@sjtu.edu.cn}{siran.li@sjtu.edu.cn}, {0000-0003-4283-273X}
    \item ZIJIU LYU: \href{mailto:zijiu.lyu@gmail.com}{zijiu.lyu@gmail.com}, 0009-0009-6456-2399
    \item HAO NI: \href{mailto:h.ni@ucl.ac.uk}{h.ni@ucl.ac.uk}, 0000-0001-5485-4376
    \item JIAJIE TAO: \href{mailto:jiajie.tao.21@ucl.ac.uk}{jiajie.tao.21@ucl.ac.uk}
\end{itemize}}
}
\thanks{$\ast$ Corresponding author.}

\date{\today}

\begin{document}

\begin{abstract}

A central question in rough path theory is characterising the law of stochastic processes on path spaces. It is established in [I. Chevyrev $\&$ T. Lyons, Characteristic functions of measures on geometric rough paths, \textit{Ann. Probab.} \textbf{44} (2016), 4049--4082] that the characteristic function of a probability measure on \emph{group-like elements}, which is a subspace of the extended tensor algebra $\tens = \prod_{n=0}^\infty\left(\R^d\right)^{\otimes n}$, uniquely determines the measure. In this work, we show that the characteristic function restricted to  special orthogonal Lie algebra $\so$(n) is sufficient to achieve this goal. The key to our arguments is an explicit algorithm --- as opposed to the non-constructive approach in [I. Chevyrev $\&$ T. Lyons, \emph{op. cit.}] --- for determining a generic element $X \in \tens$ from its generating function when restricting its domain to a sparse subspace of real tridiagonal skew-symmetric matrices. Our only assumption is that $X$ has a non-zero ROC, which relaxes the condition of having an infinite ROC in the literature. 

As an application, we propose the restricted path characteristic function distance (RPCFD), a novel distance function for probability measures on the path space that serves as the sparse counterpart of path characteristic function distance. It has enormous advantages in dimension reduction and potential in generative modeling for synthetic time series generation, validated in this paper via hypothesis testing on fractional Brownian motions.



\end{abstract}
\maketitle

\section{Introduction}\label{sec: intro}

Characterising the law of stochastic processes is of enormous significance in stochastic analysis, statistics, and machine learning \cite{icm}. The study of \emph{faithful} and \emph{parsimonious} representations of laws on the path space, in particular, has a direct impact on algorithm design and efficiency. It has found successful applications in numerical methods of SDE simulations, hypothesis testing on stochastic processes, and generative models for synthetic time series, etc. \cite{kidger2019deep,ni2023conditional, morrill2019signature, chevyrev2022signature}.

The \emph{signature} of paths, a central concept in rough path theory, can be regarded as non-commutative monomials on the path space. It uniquely determines the underlying path of finite $p$-variations up to a negligible equivalence class where $p \geq 1$ (\emph{i.e.}, the \emph{tree-like equivalence} \cite{hl, unique}). Therefore, characterising the law of stochastic processes is equivalent to characterising probability measures on the space of path signatures. 

In this paper, we consider measures on some subspace $\mathscr{V}$ (to be specified later) larger than the space of path signatures in the extended tensor algebra $\tens := \prod_{n=0}^\infty \left(\R^d\right)^{\otimes n}$.Let $\mu$ be a Radon measure on $\mathscr{V}$. We define its \emph{characteristic function} by
\begin{equation} \label{eq: intro_char}
    \Phi_\mu (M) := \mathbb{E}_{X \sim \mu} \left[ \widetilde{M}(X) \right].
\end{equation}
Here, $M$ ranges over linear mappings from $\R^d$ to $\bigoplus_{k=1}^\infty \mathfrak{su}(k)$, \emph{i.e.}, the space of skew-Hermitian matrices, and $\widetilde{M}$ denotes the canonical extension of $M$. See \eqref{tilde M, def} for its precise definition. Such a generalisation of the characteristic function appeared in \cite{cl}, and it coincides with the \emph{path characteristic function} PCF, namely the expectation of the unitary development of the stochastic process, when the measure $\mu$ is restricted to the space of path signatures \cite{ours}. 


    The key result of the seminal paper \cite{cl} by Chevyrev--Lyons ascertains that PCF determines the law of the signature of the underlying stochastic process, and PCF is determined by the expected signature provided that the latter has an \emph{infinite radius of convergence} (ROC). The proof in \cite{cl} is non-constructive and leaves open two important questions: (1), whether the entire PCF, rather than some of its simpler variants, is necessary in order to determine the law of signature, and (2), how the PCF or its variants determine the expected signature, if the expected signature exists and satisfies certain regularity conditions. Our work offers a novel and constructive approach to address these questions. We show that PCF, when restricted to $\bigoplus_{k=1}^\infty \mathfrak{so}(k)$ (\emph{i.e.}, the space of real skew-symmetric matrices), is sufficient to determine the law of signature. For the case where the expected signature has the infinite radius of convergence, our approach provides an explicit way to recover the expected signature, which can not be achieved by the approach of \cite{cl}.

To achieve these results, a key problem is to show that the family of matrix-valued functions
\begin{equation}\label{A, def, new}
\mathcal{A} := \left\{\widetilde{M}, \text{ canonical extension of } M \in \bigoplus_{k=1}^\infty \lin\left(\R^d; \mathfrak{su}(k)\right)\right\}    
\end{equation}
separates points over $\mathscr{V}\subset\tens$, in view of the Stone--Weierstrass theorem. Chevyrev--Lyons established in \cite[Theorem 4.8]{cl} that $\mathcal{A}$ separates points over the subspace of tensors with \emph{infinite radii of convergence} (ROC) by exploiting polynomial identities in the numerical linear algebra literature \cite{alg1, alg2}. Our Main Theorem~\ref{thm: main} strengthens this result in three aspects: 
\begin{enumerate}
    \item We show that $\mathcal{A}$ indeed separates points over the subspace of tensors with \emph{non-zero ROC}.
    \item We show that a (much) smaller subset $\mathcal{A}_{\mathfrak{t}} \subset \mathcal{A}$ is sufficient to achieve this goal.
    \item We give an explicit algorithm for recovering $X \in \tens$ from $\left\{\widetilde{M}(X):\,\tilm \in \mathcal{A}_{\mathfrak{t}}\right\}$, provided that $X$ has a \emph{non-zero ROC}. 
\end{enumerate}
Here and hereafter, we set \[
\mathcal{A}_{\mathfrak{t}} := \left\{\widetilde{M}, \text{ canonical extension of } M \in \bigoplus_{k=1}^\infty \lin\left(\R^d; \mathfrak{t}(k,\R)\right)\right\},
\]
where $\mathfrak{t}(k,\R)$ is the space of ``tridiagonal''  skew-symmetric real matrices:
\begin{equation} \label{t, def}
    \mathfrak{t}(k,\R) := \bigg\{M \in \R^{k \times k}:\, \text{$M+M^\top=0$; the $(i,j)$-entry of $M=0$ whenever $|i-j|\neq 1$} \bigg\}.
\end{equation}

Let us first introduce the notions of signature, development, generating function, and radius of convergence in the subsections below.

\subsection{Path signature}

Consider a path $\gamma: [0,T] \to \R^d$ of bounded $p$-variations where $1 \leq p < 2$. The \emph{signature of} $\gamma$, first introduced by K.-T. Chen in 1950's (\cite{c1, c2}), is defined by a sequence of iterated integrals:
\begin{align}
    S(\gamma)_{[0,T]} & := \left(1, \gamma^{(1)}_{[0,T]}, \ldots,\gamma^{(n)}_{[0,T]},\ldots \right), \\
    \gamma^{(n)}_{[0,T]} & := \int_{0<t_1<t_2<\ldots<t_n<T} \dd \gamma_{t_1} \otimes \cdots \otimes \dd \gamma_{t_n}\,\in \left(\R^d\right)^{\otimes n}, \label{signature, nth level}
\end{align}	
that takes values in the \emph{extended tensor algebra}:
\begin{equation}
    \tens := \prod_{n=0}^\infty\left(\R^d\right)^{\otimes n};
\end{equation}
set $\left(\R^d\right)^{\otimes 0} \equiv \R$. The iterated integral in \eqref{signature, nth level} is understood in Young's sense \cite{book}. The signature is time-reparametrisation-invariant; hence, we drop the subscript ${}_{[0,T]}$ in $S(\gamma)$ from now on.

The signature transform $\gamma \mapsto S(\gamma)$ defines a group homomorphism from the path group (with path-concatenation $\ast$ and path-inverse $\overleftarrow{\cdot}$) into $\tens$, such that
\begin{equation}\label{chen's id, group hom}
    S(\gamma\ast \eta) = S(\gamma) \otimes S(\eta) \qquad \text{and}\qquad S(\gamma) \otimes S\left(\overleftarrow{\gamma}\right) = {\bf 1} := (1,0,0,\ldots). 
\end{equation}
The first identity in \eqref{chen's id, group hom}, first established in \cite{c3}, is known as Chen's identity.

The signature contains all essential information of a path. More precisely, Chen showed in \cite{c3} that two irreducible piecewise $C^1$-paths cannot have the same signature unless they differ only by a translation or reparametrisation. This result was extended to BV-paths by Hambly--Lyons in \cite{hl}, proving that the signature uniquely determines BV-paths modulo tree-like equivalence. And we refer to \cite{unique} for the tree-like equivalence of finite $p$-variation paths where $1 < p < 2$.

Besides its theoretical importance, as a mathematically principled feature that enjoys desirable analytic-geometric properties (\emph{e.g.},  \emph{universality} and \emph{characteristicity}), the signature establishes itself as a powerful tool in machine learning applications, especially for time series data modelling. It serves as an efficient feature extractor of time series and streaming data and, when integrated with appropriate machine learning techniques, signature-based approaches have shown success in quantitative finance, early detection of disease, and online character recognition, etc. \cite{kidger2019deep,ni2023conditional, morrill2019signature, chevyrev2022signature}. See also  Lyons' ICM talk \cite{icm} and the recent works \cite{ck, lln, ours}.

\subsection{Path development}

An important tool for the study of signature is the \emph{(Cartan) development} of a path. Consider a Lie group $G$ with Lie algebra $\g$. The development of a path $\gamma:[0,T]\to\R^d$ into $G$ is a map
\begin{equation}
    \dev_\gamma: \lin\left(\R^d;\g\right) \longrightarrow G
\end{equation}
given (formally) by
\begin{equation}\label{dev, def}
\dev_\gamma(M) := I_k+\sum_{n=1}^\infty  \int_{0<t_1<t_2<\ldots<t_n<T} M\left(\dd \gamma_{t_1}\right)\bullet  \cdots \bullet M\left(\dd \gamma_{t_n}\right), 
\end{equation}
where $I_k$ is the identity of $G$.
This is well defined when $\gamma$ is a path of bounded $p$-variation for $1 \leq p < 2$. The symbol $\bullet$ denotes matrix multiplication, and $\lin\left(\R^d;\g\right)$ is the space of linear transforms from $\R^d$ to $\g$. The formula~\eqref{dev, def} is closely related to the signature of $\gamma$ defined in \eqref{signature, nth level}. 

The path development plays a central role in the study of signature inversion \cite{lx} and the uniqueness of expected signature; \emph{cf.} \emph{e.g.}, \cite{cl, hl, bg}. Chevyrev--Lyons \cite{cl} proved that, for suitably
chosen $\g$, the development is a universal, characteristic feature. Moreover, the development takes values in Lie groups of finite dimensions (independent of the path dimension $d$), thus serving as a promising tool in machine learning applications for overcoming the ``curse of dimensionality'' inherent to the path signature. See \cite{ours} and the references cited therein.

For $G = U(k)=$ the unitary group and $\g=\su(k)=$ the unitary algebra (\emph{i.e.}, the space of skew-Hermitian matrices $M$ with $M+M^*=0$), we call $\dev_\gamma$ the \emph{unitary development} of $\gamma$. In machine learning papers (\emph{e.g.}, \cite{ours, ours'}), it is also known as the \emph{unitary feature}.

\begin{remark} \label{rmk:CDE}
    Let $1 \leq p < 2$. The development of a finite $p$-variation path $\gamma: [0, T] \rightarrow \R^d$ under $M \in \lin\left(\R^d;\g\right)$ can be equivalently defined by the unique solution of the following differential equation
    \begin{align} \label{eq:CDE}
        \dd Y_t = Y_t \bullet M(\dd \gamma_t), \quad Y_0 = I_k,
    \end{align}
    where the integral is understood in Young's sense. Notice that \eqref{dev, def} coincides with the Picard iteration of \eqref{eq:CDE}.
\end{remark}

\subsection{Generating function of a tensor} \label{sec:generating_func}

Let $x\equiv \{x_n\}_{n=0}^\infty$ be a real sequence. Its \emph{generating function} is defined as the (formal) complex series:
\begin{align}\label{gen fn for real seq}
    \phi_x: \C \longrightarrow \C, \qquad \phi_x(\lambda) := \sum_{n=0}^\infty x_n \lambda^n.
\end{align}
For $\lambda \in \C$ with $|\lambda|<\rho(x):=$ the radius of convergence of the above power series, we may differentiate the series term-by-term. Thus, with $\phi_x$ given, we may recover $x$ via
\begin{equation}\label{recover x, R case}
    x_n = \frac{1}{n!} \frac{\dd^n}{\dd \lambda^n}\phi_x(\lambda)\bigg|_{\lambda = 0}.
\end{equation}
Now let us generalise the above construction to tensors.
\begin{definition}\label{def, gen fn}
    For each given $k= 1,2,3,\ldots$ and tensor $X \in \tens$, the \underline{generating function} \underline{of $X$ of order $k$} is defined by the mapping:
\begin{align} \label{gen fn, formula}
    \Phi_X^{(k)}: \lin\left(\R^d; \su(k)\right) \,\longrightarrow\, \matk, \qquad \Phi_X^{(k)}(M) := \tilm(X).
\end{align}
Here $\matk \cong \C^{k\times k}$ is the algebra of $k\times k$ complex matrices, and $\tilm: \tens \to \matk$ is the \underline{canonical extension of $M$} defined as follows:
\begin{multline} \label{tilde M, def}
\text{$\tilm$ is the algebra homomorphism determined by $\tilm(\mathbf{1}) = I_k$ and that } \\
\tilm\left(e_{i_1} \otimes  \cdots \otimes e_{i_\ell}\right):= M\left(e_{i_1}\right) \bullet \cdots \bullet M\left(e_{i_\ell}\right) \text{ for any } \ell \in \mathbb{N} \text{ and } i_1, \ldots, i_\ell \in \{1,\ldots,d\}.
\end{multline}
The \underline{generating function of $X$} is defined by:
\begin{eqnarray}
&&\Phi_X: \bigoplus_{k=1}^\infty \lin\left(\R^d; \su(k)\right) \,\longrightarrow\, \bigoplus_{k=1}^\infty\matk,\\
&& \Phi_X\Big|_{\lin\left(\R^d; \su(k)\right)} := \Phi_{X}^{(k)},
\end{eqnarray}
which respects the natural grading in $k \geq 1$.
\end{definition}

Notice that the defining equation~\eqref{gen fn, formula} for the generating function is formal --- consider the projection of a tensor onto the $n^{\text{th}}$ level: 
\begin{equation}\label{proj}
\pi_n: \tens \longrightarrow \left(\R^d\right)^{\otimes n}.
\end{equation}
In view of \eqref{tilde M, def}, the action of $\tilm$ on $\pi_n(X)$ involves a product of $n$ matrices. To make \eqref{gen fn, formula} rigorous, we require that each $n$-fold product in \eqref{tilde M, def} converges in $\matk$ as the tensor degree $n\to\infty$. To fix the idea,  in this paper we always adopt the Hilbert--Schmidt norm on $\matk$.

Definition~\ref{def, gen fn} of the generating function can be extended to any Lie algebra $\mathfrak{g}$ and normed unital algebra $\mathfrak{a}$. One defines $\Phi_X: \mathcal{L}\left(\R^d; \mathfrak{g}\right) \to \mathfrak{a}$ again via $\Phi_X(M) := \tilm(X)$, where $\tilm$ is given by \eqref{tilde M, def} with $\bullet$ being the multiplication in $\mathfrak{a}$. If $X \in \mathcal{G}\left(\R^d\right)$, then $\Phi_X$ takes values in $G$, the Lie group of $\mathfrak{g}$. On the other hand, if $X$ is in ${\rm Lie}\left[\mathcal{G}\left(\R^d\right)\right]$, the Lie algebra of $\mathcal{G}\left(\R^d\right)$, then $\Phi(X) \in \mathfrak{g}$.\footnote{The Lie algebra of $\mathcal{G}\left(\R^d\right)$ is denoted as $\mathcal{L}((\R^d))$ in \cite{book}. We avoid using this notation and reserve  $\mathcal{L}$ for the space of linear maps.} Here and hereafter,
\begin{equation*}
\mathcal{G}\left(\R^d\right) := \left\{\text{group-like elements in $\tens$}\right\},  
\end{equation*}
and ${\rm Lie}\left[\mathcal{G}\left(\R^d\right)\right]$ is spanned by $\left[e_{i_1}, \left[e_{i_2}, \cdots, \left[e_{i_{\ell-1}}, e_{i_{\ell}}\right]\right]\right]$ where $\{e_1, \ldots, e_d\}$ is a basis for $\R^d$. See \cite[\S 2.2.4]{book}.

\begin{example}\label{ex: 1d}
When $d=1$ and $k=1$, the 1-dimensional matrix algebra $\su(1)$ equals $\2\R$. We view $\R^{\otimes n}\cong \R$ and hence the tensor algebra satisfies 
$$
\T(\R) \cong \prod_{n=0}^\infty \R = \text{the space of real sequences}.
$$ 
Thus, for $x\in \T(\R)$ one identifies $\pi_n(x)=x_n 1^{\otimes n} \in \R^{\otimes n}$ and 
\begin{equation}\label{series, notation}
x=(x_0,x_1, x_2,\ldots) \sim \sum_{n=0}^\infty x_n 1^{\otimes n}.
\end{equation}
In this way, for each $\lambda \in \lin(\R; \2\R)$ we view $$\widetilde{\lambda}\left(x_n 1^{\otimes n} \right) = x_n \overbrace{\lambda(1)\cdots \lambda(1)}^{n \text{ times }} \equiv x_n \lambda^n,$$ with the obvious abuse of notation $\lambda(1) \equiv \lambda$. 

Therefore, in view of the series notation in \eqref{series, notation}, we may express $\phi_x(\lambda) = \sum_{n=0}^\infty x_n \lambda^n=\Phi_x(\lambda)$. That is, the generating function for tensors over $\R^1$ into 1-dimensional unitary matrix algebra agrees with the classical notation in \eqref{gen fn for real seq} for sequences. This also explains why we need to consider \emph{complex} series for  \eqref{gen fn for real seq}.  
\end{example}

Observe the structural resemblance of the integrand of the development in \eqref{dev, def} with the canonical extension in \eqref{tilde M, def} of the generating function. Indeed, we have the following result, in which $\su(k)$ can be replaced by any finite-dimensional Lie algebra.
\begin{proposition}\label{propn: dev and sig}
Let $1 \leq p < 2$. The development and signature of a finite $p$-variation path $\gamma:[0,T]\to\R^d$ satisfies 
\begin{equation*}
\dev_\gamma(M) = \tilm \big(S(\gamma)\big)
\end{equation*}
for any $k=1,2,3,\ldots$ and $M \in \lin\left(\R^d; \su(k)\right)$.
\end{proposition}

In passing, we remark that Definition~\ref{def, gen fn} for the generating function is closely related to the notion of the \emph{characteristic function} $\phi_X$ for a $G$-valued random variable $X$, where $G$ is a topological group (see Heyer \cite{heyer} and Chevyrev--Lyons \cite[\S 1]{cl}):
\begin{equation}\label{char fn, def}
\phi_X (M):= \mathbb{E}\left[\tilm(X)\right]\qquad \text{where $M$ is any unitary representation of $G$.}
\end{equation}
In light of Proposition~\ref{propn: dev and sig}, if $X=S(\gamma)$ for some random finite $p$-variation path $\gamma$ where $1 \leq p < 2$, then $\phi_X(M)$ equals the expectation of the development $\dev_\gamma(M)$. Here $G$ can be taken as $\mathcal{G}\left(\R^d\right):=\left\{\text{group-like elements in $\tens$}\right\}.$ See \cite[\S 3]{cl}.

\subsection{Radius of convergence of a tensor}\label{sec: ROC}

One central question concerning the generating function of a tensor is whether the analogue of \eqref{recover x, R case} holds. 

\begin{question}\label{q: mom prob}
Given $\Phi_X:\bigoplus_{k=1}^\infty \lin\left(\R^d; \su(k)\right) \to \bigoplus_{k=1}^\infty \matk$ the generating function of $X$, graded by $k \geq 1$. Can we solve for $X$?
\end{question} 
\noindent
Of particular interests are the following choices of $X$:
\begin{itemize}
\item
Signature --- $X=S(\gamma)$ for some path $\gamma$, \emph{e.g.}, of bounded $p$-variation for $1\leq p <2$;
\item
Expected signature --- $X = \mathbb{E}[S(\gamma)]$ for some \emph{random} path $\gamma$.
\end{itemize}
By virtue of Proposition~\ref{propn: dev and sig}, in the above cases Question~\ref{q: mom prob} specialises to the following:
\begin{question}\label{q: sig/expsig}
Does the (expected, resp.) development uniquely determines the (expected, resp.) signature of a path?
\end{question}
\noindent
Notice that certain conditions must be imposed on $X$ to get an affirmative answer to Question~\ref{q: mom prob}. Even in the simplest case of real sequences (Example~\ref{ex: 1d}), to solve for each $x_n$ via \eqref{recover x, R case}, the radius of convergence of the power series in \eqref{gen fn for real seq} must be non-zero. 

The notion of \emph{radius of convergence} can be extended to the elements of $\tens$. This shall play a central role in the investigation of Questions~\ref{q: mom prob} $\&$ \ref{q: sig/expsig}. Throughout this paper, to fix the ideas, we equip $\R^d$ with the $\ell_1$-norm, namely that $\|a\|:=\sum_{i=1}^d |a_i|$ for $a=(a_1,\ldots,a_d)^\top$, and equip $\left(\R^d\right)^{\otimes n}$ with the projective tensor norm for each $n=1,2,3,\ldots$, denoted by $\|\cdot\|$ again\footnote{See \cite[pp.136--139]{hl} and \cite[p.3]{nixu} for discussions on these norms. However, one may also define $\E_\infty$ as the completion of $\bigoplus_{n=0}^\infty \left(\R^d\right)^{\otimes n}$ under scalar
dilations of any system of sub-multiplicative norms, as per \cite[Remark~5.3]{cl}.}. Also recall the level-$n$ projection operator $\pi_n$ from \eqref{proj}. Then, we introduce:
\begin{definition}\label{def: ROC}
The \underline{radius of convergence} (ROC) of a tensor $X \in \tens$, written as $\rho(X)$, is defined by the radius of convergence of the (complex) power series 
\begin{align*}
\lambda \longmapsto \sum_{n=0}^\infty \left\|\pi_n(X)\right\| \lambda^n.
\end{align*}
Note that $\rho(X) \in [0,\infty]$. For any $R\geq 0$, we further set
\begin{equation*}
    \E_R:= \left\{X \in \tens:\, \sum_{n=0}^\infty \left\|\pi_n(X)\right\| R^n \text{ is convergent}\right\},\quad \hatE:= \bigcup_{R>0}\E_R,\quad  \E_\infty := \bigcap_{R>0} \E_R.
\end{equation*}
\end{definition}

Note that if $X \in \E_R$, then $\rho(X) \geq R$. Also, $\E_R$ is a Polish space under the metric $\sum_{n=0}^\infty \left\|\pi_n(X)-\pi_n(Y)\right\| R^n$.

The radius of convergence $\rho(X)$ of a tensor $X$ is closely related to its generating function $\Phi_X$ defined in Definition~\ref{def, gen fn}, as shown by the following proposition.

\begin{proposition}\label{propn: ROC} Given $X \in \tens$ and $k =1,2,3,\ldots$.
\begin{enumerate}
\item
If $X \in \E_\infty$, \emph{i.e.}, $\rho(X)=\infty$, then the generating function $\Phi_X$ in \eqref{gen fn, formula} is well defined on the domain $\lin\left(\R^d; \su(k)\right)$.
\item
If $X \in \E_R$,  then  $\Phi_X$ is well defined on  $\left\{M \in \lin\left(\R^d; \su(k)\right):\, \text{the operator norm of $M$} < R\right\}$.
\item
$X\in \E_\infty$ if and only if it lies in the closure of $\bigoplus_{n=0}^\infty \left(\R^d\right)^{\otimes n}$ with respect to the coarsest topology for which the following holds: for any normed algebra $\mathcal{A}$ and any algebra homomorphism $M \in {\rm Hom}\left(\R^d, \mathcal{A}\right)$,  the canonical extension $\tilm: \tens \to \mathcal{A}$ is continuous.
\end{enumerate}
\end{proposition}

\begin{proof}
Clearly $(2) \Rightarrow (1)$. To see (2), we simply note that $\left\|\widetilde{M}(\pi_n(X))\right\| \leq \|\pi_n(X)\| \cdot \|M\|_{\rm op}^n$ for each $n \in \mathbb{N}$, where $\|M\|_{\rm op}$ is the operator norm of $M$. Finally, (3) follows from \cite[Corollary~2.5]{cl}, with $\Psi$ chosen to be the singleton of the $\ell_1$-norm on $\R^d$ therein.  \end{proof}

\subsection{Main Theorem}

Our main result of this paper provides an affirmative answer to Question~\ref{q: mom prob}, under the assumption $\rho(X)>0$ (equivalently, $X \in \hatE$). In fact, we establish an explicit algorithm to solve for $X$ from its generating function $\Phi_X$:

\begin{theorem}\label{thm: main}
Given the generating function (graded by $k \geq 1$)
$$
\Phi_X : \bigoplus_{k=1}^\infty \lin\left(\R^d;\su(k)\right) \longrightarrow \bigoplus_{k=1}^\infty \matk
$$ 
of some $X\in \hatE$, \emph{i.e.}, a tensor with a non-zero radius of convergence. Then, one can recover $X$ explicitly by the following procedure:

(i) Write $$X = \sum X^W e_{w_1} \otimes \cdots \otimes e_{w_n},\qquad W \equiv \left(w_1, \ldots, w_n\right),$$ where $\{e_1, \ldots, e_d\}$ is the Cartesian basis for $\R^d$ and the summation is taken over all possible $n \in \mathbb{N}$ and ordered $n$-tuples $W \in \dset^n$. It suffices to determine $X^W$ for each $W$.

(ii) To this end, fix any $k \geq n+1$ and take the linear mappings $$\mw_1, \ldots, \mw_d \in \lin\left(\R^d; \mathfrak{t}(k,\R)\right)$$ defined for each $i,j \in \{1,\ldots,d\}$ as follows:
\begin{equation*}
\mw_i(e_j) := 
\begin{cases}
\sum'_{t} \left(E^{p^W(i|t)}_ {p^W(i|t)+1} - E^{p^W(i|t)+1}_{p^W(i|t)}\right) \quad \text{ if } i = j \in \{w_1, \ldots, w_n\},\\
0\qquad \text{ if } i \notin \{w_1, \ldots, w_n\} \text{ or } i \neq j.
\end{cases}
\end{equation*}
Here, $p^W(i|t)$ is the location where $i$ appears for the $t^{\text{th}}$ time in $W=(w_1, \ldots, w_n)$, and $\sum'_{t}$ means summation over $t$ from $1$ up to the total number of times that $i$ appears in $W$. The symbol $E^a_b$ designates the $k\times k$ matrix whose $(a,b)$-entry is $1$ and all the other entries are $0$.

(iii) Form the linear combination 
\begin{equation*}
\mwt = \sum_{i=1}^d \theta_i \mw_i,
\end{equation*}
where $\theta = (\theta_1, \ldots, \theta_d)^\top \in \R^d$ are ``parameters''.  Then
\begin{equation}\label{conclusion, in theorem}
X^W = \frac{1}{\cw} \times \left\{\text{the $(1,n+1)$-entry of the matrix  } \frac{\p^n\Phi \left(\mwt\right)}{\p\theta_{w_1}\cdots\p\theta_{w_n}}\Bigg|_{\theta=0}\right\}, 
\end{equation}
where the combinatorial constant $\cw$ equals the product over each $i \in \dset$ of the factorial of the number of times that $i$ appears in $W$.
\end{theorem}

Despite the lengthiness of the statement of Theorem~\ref{thm: main}, we emphasise its explicit and constructive nature.  Every input (\emph{i.e.},  constant $\cw$ and matrices $\mwt$) in the formula~\eqref{conclusion, in theorem} depends only on the word $W$ and is purely combinatorial. Detailed explanations for the notations can be found in \S\ref{sec: notation} below.

An immediate consequence of Theorem~\ref{thm: main} is the \emph{characteristic property} of the generating function. That is, generating functions (viewed as the map $X \mapsto \Phi_X$) separate points in $\hatE$.
\begin{corollary}\label{cor: char}
    For any $X \in \hatE$ with $\pi_n(X) \neq 0$ and any $k \geq n+1$, one can find $M \in \lin\left(\R^d; \mathfrak{t}(k,\R)\right)$ such that $\Phi_X(M) \neq 0$. Here $\Phi_X$ is the generating function of $X$ as in Definition~\ref{def, gen fn}, and $\mathfrak{t}(k,\R)$ is the space of tridiagonal skew-symmetric $k\times k$ real matrices as in \eqref{t, def}.
\end{corollary}

This result should be compared with Chevyrev--Lyons \cite[Theorem~4.8]{cl}, which ascertains that for any $X \in \E^\infty$ with $\pi_n(X) \neq 0$, one can find $M \in \lin\left(\R^d; \mathfrak{sp}(m)\right)$ for some $m \geq \max \left\{2, \frac{n}{3}\right\}$ such that $\Phi_X(M) \neq 0$. Here $\mathfrak{sp}(m)$ is the Lie algebra of (the real form of the) symplectic group of $2m \times 2m$ real matrices, which contains $\mathfrak{t}(2m,\R)$ as a subset. The proof in \cite{cl}, based on polynomial identities in the numerical linear algebra literature \cite{alg1, alg2}, is non-constructive by nature.

Given the above discussions on radius of convergence, we can say more about Question~\ref{q: sig/expsig}. 

\begin{itemize}
    \item 
    As shown in Chevyrev--Lyons \cite[Corollary~2.5 and \S 5]{cl}, the signature of a $p$-rough path for any $p\geq 1$ has an infinite radius of convergence.  Hence, its moment generating function uniquely determines the signature. In particular, for $p \in [1, 2[$ the signature can be defined as in \eqref{signature, nth level} in the sense of Young's integrals, and the signature of any geometric rough path in $\R^d$ lies in $\mathcal{G}\left(\R^d\right)$, the set of group-like elements in $\tens$. In this case, the uniqueness of unitary development holds.

    \item 
    The case for the expected signature is more complicated. The infinity of radius of convergence for the expected signature has been verified for several popular stochastic processes on a fixed time horizon, \emph{e.g.}, fractional Brownian motions with Hurst parameter in $\left[\frac{1}{2},1\right[$ (\cite{hurst}), as well as certain Gaussian and Markovian rough paths (\cite[Examples~6.7 and 6.8]{cl}). 
    It is, nonetheless, considerably challenging to study processes up to a \emph{random} time. The only known results in this regard, to the best of our knowledge, pertain to the Brownian motion up to the first exit time from a bounded domain $\Omega \subset \R^d$ (\emph{a.k.a.} the ``stopped Brownian motion''). Using PDE approaches pioneered by Lyons--Ni \cite{phd}, it has been established in \cite{bdmn, lini} that the expected signature of stopped Brownian motions on $C^{2,\alpha}$-domains in $\R^d$ for $2 \leq d \leq 8$ has a \emph{finite} radius of convergence. 
\end{itemize}

\subsection{Applications to numerical analysis and machine learning}

The use of signature and path development in machine learning is on the rise, due to the characteristic and universal properties. The signature serves as an efficient feature extractor of time series and streaming data. When integrated with appropriate machine learning techniques (\emph{e.g.}, deep learning and gradient boosting trees), signature-based approaches have shown success across various fields --- quantitative finance, early detection of disease, and online character recognition, etc. \cite{kidger2019deep,ni2023conditional, morrill2019signature, chevyrev2022signature}. 

However, the signature is prone to \emph{the curse of dimensionality}, especially in tasks involving high-dimensional time series. To this end, \cite{ours'} proposed the \emph{path development layer}, which is data-adaptive and significantly reduces dimensions while preserving essential information of the signature feature. Built on top of the expected unitary development, \cite{ours} further proposed the so-called \emph{path characteristic function (PCF)}, which proves to be a general metric on the space of distributions on the unparameterised path space. Successful applications of PCF have been found in hypothesis testing on time series distributions and in synthetic time series generation.

The theoretical justification of the PCF distance is essentially based upon the characteristic property of the (expected) unitary development. In Theorem~\ref{thm: main} and its proof, we specify a reduced set of sparse linear maps in $\mathcal{L}\left(\R^d; \g\right)$, under which the corresponding unitary development uniquely determines the unparameterised path. It motivates us to consider the ``restricted PCF distance'', which characterises the law of time series while reducing the computation complexity compared with the PCF distance in \cite{ours}. This will be introduced in \S\ref{sec: ht}, with demonstration of its efficacy in hypothesis testing on fractional Brownian motion.

\subsection{Organisation}

The remaining parts of the paper are organised as follows.

In \S\ref{sec: notation}, we introduce some notations that shall be used throughout this paper. Then, in \S\ref{sec: der}, a combinatorial identity on derivatives of the unitary development of special matrices is presented. Our Main Theorem~\ref{thm: main} is proved in \S\ref{sec: proof}. The following section \S\ref{sec: ESig} is devoted to the theoretical implications of Main Theorem~\ref{thm: main}, \emph{i.e.}, the determination of measures on subspaces of $\tens$ and path spaces. Finally, in \S\ref{sec: ht}, we investigate the application of Main Theorem~\ref{thm: main} to hypothesis testing, via the concept of restricted path characteristic function distance introduced therein. Appendix~\ref{details} presents the details of the numerical experiments on hypothesis testing.

\section{Notations and Nomenclature}\label{sec: notation}

In this preliminary section, we introduce some notations that will be used throughout. We assume that all paths in consideration take values in $\R^d$, with $d$ fixed once and for all. 

\subsection{Words}
We denote
\begin{align*}
\wordn := \Big\{W:\, \text{$W$ is a word of length $n$ from the alphabet $\{1,\ldots,d\}$}\Big\}.
\end{align*}
That is, $W \in \wordn$ is an ordered $n$-tuple $W = \left(w_{1}, \ldots, w_{n}\right)$ where $w_{i} \in \dset$ for each $i \in \nset$. As a set $\wordn = \dset^n$; each $w_i$ appearing in $W$ is called a \emph{letter} of $W$.

A generic tensor $X \in \tens$ can be expressed as the (formal) sum:
\begin{align}\label{tensor sum}
X = \sum_{n=0}^\infty \,\, \sum_{W \in \wordn} X^W e_W
\end{align}
where, for each $W = \left(w_{1}, \ldots, w_{n}\right)$ as above, $X^W \in \R$ is a real coefficient and $e_W \equiv e_{w_1} \otimes \cdots \otimes e_{w_n}$. 

From an alternative perspective, we can view each $W \in \wordn$ as a \emph{multi-index} with valency $|W|=n$. One may take partial derivatives via 
\begin{align*}
\frac{\p^n}{\p x^W} = \frac{\p^n}{\p x^{w_1} \cdots \p x^{w_n}},
\end{align*}
with possible repetitions in $\{w_1, \ldots, w_n\}$.

Denote by $\Sigma_n$ the permutation group of $n$ letters, \emph{i.e.}, the set of bijections from $\nset$ to itself. One has the natural group action $$\Sigma_n \times \wordn \rightarrow \wordn$$ given by
\begin{align*}
\sigma \cdot \left(w_1, \ldots, w_n\right) := \left(w_{\sigma(1)}, \ldots, w_{\sigma(n)}\right).
\end{align*}
One may define an equivalence relation on $\wordn$:
\begin{align*}
I \sim J \qquad \text{if and only if}\qquad I = \sigma(J) \text{ for some } \sigma \in \Sigma_n.
\end{align*}
That is, $I \sim J$ means that they differ by a permutation, hence have the same length and same (unordered) indices/letters.

\subsection{Multiplicity and position functions}
For each word $W \in \wordn$ and each letter of it, we define the multiplicity and position functions:
\begin{definition}\label{def: r and p}
Given $W \in \wordn$ where $n =1,2,3,\ldots$. 
\begin{enumerate}
\item
For $i \in \dset$, we set $r^W(i) :=$  the multiplicity of $i$ in $W$, \emph{i.e.}, the number of times that $i$ appears in $W$.
\item
For  $i \in \dset$ and $t \in \left\{1,\ldots,r^W(i)\right\}$, we set $p^{W}(i|t):=$ the position that $i$ appears in $W$ for the $t^{\text{th}}$ time. We put $p^{W}(i|t)=\emptyset$ if $r^W(i)=0$.
\end{enumerate} 
\end{definition}

\begin{example}\label{example 12221}
Consider the case $d=4$, $n=5$, and $W=(12221)$. Then
\begin{align*}
r^W(1) = 2, \quad r^W(2)=3, \quad r^W(3)=r^W(4)=0;
\end{align*}
as well as
\begin{align*}
&p^{W}(1|1)=1, \quad p^W(1|2)=5,
\quad p^{W}(2|1)=2,\quad p^{W}(2|2)=3,\quad p^{W}(2|3)=4,\\
&\qquad \text{and } \quad p^W(i|t) = \emptyset \quad \text{ for all the other $i \in \{1,2,3,4\}$ and $t \in \left\{1,\ldots,r^W(i)\right\}$}.
\end{align*}
\end{example}

Note that for any $W \in \wordn$,  $i \in \dset$, and $t \in \left\{1,\ldots,r^W(i)\right\}$, we have
\begin{align*}
r^W(i) \in \{0,1,\ldots, n\} \quad \text{and} \quad p^{W}(i|t) \in \nset. 
\end{align*}
The function $t\mapsto p^{W}(i|t)$ is strictly increasing for fixed $i$. One also has the simple identity:
\begin{equation*}
\sum_{i=1}^n r^W(i) = n \qquad \text{ for any } W \in \wordn. 
\end{equation*}

\subsection{Matrices}
Fix $k \in \mathbb{N}$ and assume that all the matrices in this paragraph are of size $k \times k$. We use the symbol $E^a_{b}$ ($a,b \in \kset$) to denote the \emph{elementary matrices}:
\begin{align*}
E^a_{b} \in \matk \quad \text{such that the $(a',b')$-entry is $\delta_a^{a'}\delta^b_{b'}$ for any $a',b' \in \kset$}. 
\end{align*}
That is, $E^a_{b}$ is the matrix whose $(a,b)$-entry is $1$ and all the other entries are $0$.

The following identity will be used repetitively:
\begin{equation}\label{id on E}
E^a_{b} \bullet E^p_{q} = \delta_b^{p} E^a_{q}\qquad\text{for any } a,b,p,q \in \kset.
\end{equation}
Throughout we use  $\bullet$ to denote the matrix multiplication. Notice also that 
\begin{align*}
\text{$\left\{E^a_{b}-E^b_{a}:\, 1\leq a < b \leq k\right\}$ is a natural basis for $\so(k,\R)$.}
\end{align*}
Here $\so(k,\R) \subset \su(k)$ is the Lie algebra of real skew-symmetric $k\times k$ matrices, which has dimension $\frac{k(k-1)}{2}$.
The matrices $\mw_1, \ldots, \mw_d$ in Theorem~\ref{thm: main} are central to our constructions. Here they are chosen to take values in $\so(k,\R)$, but see also Remark~\ref{inv method remark}.

\section{Derivatives of the development with respect to hyperparameters}\label{sec: der}

Here and hereafter, we take a tensor $X \in \hatE$.

Define for a hyperparameter $\theta = (\theta_1, \ldots, \theta_d) \in \R^d$ the matrix
\begin{align}\label{M theta expansion}
\mt := \sum_{i=1}^d \theta_i M_i,
\end{align}
where $M_1, \ldots, M_d$ are arbitrary linear maps in $\lin\left(\R^d;\su(k)\right)$. The generating function of $X$ acting on $\mt$ can be computed as follows:
\begin{align}\label{dev expression, long}
\Phi_X (\mt) &:= \widetilde{\sum_{i=1}^d \theta_i M_i}\, \left(\sum_{n=0}^\infty \,\, \sum_{W \in \wordn} X^W e_W\right)\nonumber\\
&=\sum_{n=0}^\infty \,\, \sum_{W = \left(w_{1}, \ldots, w_{n}\right) \in \wordn} X^W \left\{\sum_{i_1=1}^d\theta_{i_1}M_{i_1}\left(e_{w_1}\right)\right\}\bullet \cdots \bullet \left\{\sum_{i_n=1}^d\theta_{i_n}M_{i_n}\left(e_{w_n}\right)\right\}\nonumber\\
&=\sum_{n=0}^\infty \,\, \sum_{W = \left(w_{1}, \ldots, w_{n}\right) \in \wordn} X^W \cdot\nonumber\\
&\qquad\qquad\qquad \cdot \left\{\sum_{I=\left(i_1, \ldots, i_n\right) \in  \wordn}\Big[\theta_{i_1}\cdots \theta_{i_n}\Big] M_{i_1}\left(e_{w_1}\right)\bullet \cdots \bullet M_{i_n}\left(e_{w_n}\right) \right\}.
\end{align}
The first equality follows from \eqref{tensor sum}, the second is deduced from the definition of $\tilm$ (with dummy indices $(i_1, \ldots, i_n)$ introduced), and the final one  holds by a direct rearrangement.

Now we are interested in computing the derivative of $\Phi_X (\mt)$ with respect to $\theta^J$ evaluated at $\theta = 0 \in \R^d$, where $J=\left(j_1 \ldots, j_L\right) \in \word(L)$ is a multi-index of valency $|J|=L$. That is, we shall find the matrix
\begin{equation}\label{to find}
\frac{\p^{|J|}}{\p\theta^J}\Phi_X \left(\mt\right)\bigg|_{\theta=0} \,\in\matn.
\end{equation}
 The indices $j_1,\ldots,j_L$ may have repetitions. Notice that for $X \in \hatE$, \eqref{to find} is always well defined, since $\frac{\p^{|J|}}{\p\theta^J}\Phi_X\left(\mt\right) $ is infinitely differentiable in $\theta$ for $|\theta|$ sufficiently small, \emph{i.e.}, when the operator norm of $\mt$ is less than $\rho(X)$.

To this end, notice that the last line of \eqref{dev expression, long} lies in $$\matn\otimes \R[\theta_{1}, \ldots, \theta_d];$$ \emph{i.e.}, it is a matrix-valued $d$-variate polynomial in $\theta_{1}, \ldots, \theta_d$. Thus, the nontrivial contributions to the quantity in \eqref{to find} arise only from those terms in the right-most side of \eqref{dev expression, long} with
\begin{itemize}
\item
words $W = \left(w_{1}, \ldots, w_{n}\right) $ of length $n=L$; and
\item
multi-indices $I:=(i_1, \ldots, i_n)\in \dset^L$ which differ from $J=\left(j_1 \ldots, j_L\right)$ only by a permutation; \emph{i.e.}, there exists $\sigma \in \Sigma_n$ such that $I=\sigma(J)$, namely that $j_{\sigma(\nu)}=i_\nu$ for each $\nu \in \nset$.
\end{itemize}
Observe also that for any $I=\sigma(J)$ it holds that
\begin{align*}
\frac{\p^{|J|}}{\p\theta^J} \bigg|_{\theta=0} \theta_{i_1}\cdots \theta_{i_n} = \cj,
\end{align*}
where $\cj$ is the combinatorial constant
\begin{equation}\label{const cj}
\cj := \prod_{\alpha=1}^{|J|} \left[r^J\left(j_\alpha\right)!\right] = \prod_{i=1}^d\left[r^J\left(i\right)!\right].
\end{equation}
Here $r^J(i)$ is the multiplicity of $i$ in the word $J$ as before, and the convention $0!=1$ is adopted. Note too that $\cj \equiv \mathcal{C}_I$ whenever $I\sim J$.

Summarising the above arguments, we obtain for any multi-index/word $J$  the following equality in $\matk$:
\begin{align}\label{key 1, derivative expression}
\frac{1}{\cj} \frac{\p^{|J|}\Phi_X (\mt)}{\p\theta^J}\Bigg|_{\theta=0} =\sum_{\WW  \in \word ({|J|})} X^\WW  \left\{\sum_{\{I\in \word ({|J|}):\, I \sim J \} } M_{i_1}\left(e_{\ww_1}\right)\bullet \cdots \bullet M_{i_{|J|}}\left(e_{\ww_{|J|}}\right) \right\},
\end{align}
where $\WW=\left(\ww_{1}, \ldots, \ww_{{|J|}}\right)$ and $I =\left(i_1, \ldots, i_{|J|}\right)$. We use the dummy variable $\WW$ here to avoid potential confusion from the subsequent developments. Note that $\cj$ is independent of both $\WW$ and $I$. 
The matrices $M_1, \ldots, M_d \in \lin\left(\R^d;\su(k)\right)$ in \eqref{M theta expansion}, as well as the size of the matrices $k$, are completely arbitrary in  \eqref{key 1, derivative expression}.

\section{Determination of tensor from the development} \label{sec: proof}
 
Now we are at the stage of proving our Main Theorem~\ref{thm: main}. We endeavour to explain the motivations for the choices of mathematical objects that appear in the statement of the theorem.

\begin{proof}[Proof of Theorem~\ref{thm: main}]

The arguments are divided into five steps below.  

\smallskip
\noindent
{\bf Step~1.}  Fix arbitrary $n \in \mathbb{N}$ and $W= \left(w_1, \ldots, w_n\right) \in \wordn$. We aim to solve for the coefficient $X^W \in \R$ in front of $e_W$; see \eqref{tensor sum} for the notations. 

In view of the comments at the end of \S\ref{sec: der}, we have the freedom of choosing $\{M_i\}_{1\leq i \leq d}$ in \eqref{M theta expansion} that depends on $W$ and allows us to determine $X^W$ from \eqref{key 1, derivative expression}. We shall designate such choices of $M_i$ as $\mw_i$, in order to emphasise that they depend only on $W$. For hyperparameters $$\theta = (\theta_1, \ldots, \theta_d) \in \R^d$$ we write
\begin{align*}
\mwt = \sum_{i=1}^d \theta_i \mw_i.
\end{align*}
Here and hereafter, we require the size of the matrices to satisfy $k \geq n+1$.

\smallskip
\noindent
{\bf Step~2.} Our starting point is the formula~\eqref{key 1, derivative expression} obtained in \S\ref{sec: der} above. Let us rewrite the inner summands thereof by
\begin{equation}\label{notation, braket}
\bbra I, \WW \kket_{M} :=  M_{i_1}\left(e_{\ww_1}\right)\bullet \cdots \bullet M_{i_{n}}\left(e_{\ww_{n}}\right) \in \matk
\end{equation}
for $\WW=\left(\ww_{1}, \ldots, \ww_{{n}}\right)$ and $I =\left(i_1, \ldots, i_{n}\right)$. 
We then recast \eqref{key 1, derivative expression} with $W$ in place of $J$ as follows: 
\begin{align}\label{two sums}
\frac{1}{\cw} \frac{\p^n\Phi_X \left(\mt\right)}{\p\theta^W}\Bigg|_{\theta=0} = \sum_{\left\{\left(\WW, I\right) \in \wordn \times \wordn :\, I\sim W\right\}} X^\WW \bbra I,\WW\kket_M. 
\end{align}

Notice that the comments ensuing \eqref{to find} explains the necessity of the assumption $X \in \hatE$. Also, by virtue of \eqref{key 1, derivative expression}, non-zero contributions to the summation on the right-hand side of \eqref{two sums} arises only from the pairs $(I,\WW)$ such that $I \sim W$.

\smallskip
\noindent
{\bf Step~3.} Next let us motivate our choice for $\left\{M_i=\mw_i\right\}_{1 \leq i \leq d}$. The essential property required for  $\left\{\mw_i\right\}_{1 \leq i \leq d}$ is the following: 
\begin{align}\label{key2}
&\text{The $(1,n+1)$-entry of } \bbra I, \WW\kket_{\mw} = \begin{cases}
1\qquad \text{if } I=W=\WW,
\\
0\qquad \text{if otherwise},
\end{cases}\nonumber\\
&\qquad\qquad\qquad \text{ for any $I,\WW \in \wordn$ such that $I \sim W$}. 
\end{align} 

Indeed, assuming \eqref{key2},  we find that all the terms on the right-hand side of \eqref{two sums} have trivial $(1,n+1)$-entry unless $I=W=\WW$. 
Recall the definition of $\cw$ from \eqref{const cj}; we then have
\begin{equation}\label{final conclusion}
X^W = \frac{1}{\prod_{i=1}^d\left[r^W\left(i\right)!\right]} \times \left\{\text{the $(1,n+1)$-entry of  } \frac{\p^n\Phi_X \left(\mwt\right)}{\p\theta^W}\Bigg|_{\theta=0}\right\},
\end{equation}
which is precisely the formula~\eqref{conclusion, in theorem} to be verified.

\smallskip
\noindent
{\bf Step~4.} We specify $\left\{\mw_i\right\}_{1 \leq i \leq d}$ as in the statement of the theorem, namely that
\begin{equation}\label{key3, Mi}
\mw_i(e_j) := \begin{cases}
\delta_{ij}\left\{ \sum_{t=1}^{r^W(i)} \left(E^{p^W(i|t)}_ {p^W(i|t)+1} - E^{p^W(i|t)+1}_{p^W(i|t)}\right) \right\}\quad \text{ if } r^W(i) \neq 0,\\
0\qquad \text{ if } r^W(i)=0.
\end{cases}
\end{equation}
 We reiterate that these matrices lie in $\lin\left(\R^d; \so(k,\R) \subset \su(k)\right)$, where $k \geq n+1$ is arbitrary.

Recall here that $E^a_{b} \in \matk$ are the elementary matrices. The multiplicity function $r^W(i)$ and position function $p^W(i|t)$ are as in Definition~\ref{def: r and p}. Note too that $\mw_i(e_j)$ have only non-zero entries (which are $\pm$) on the super- and sub-diagonal, \emph{i.e.}, at $(q, q\pm 1)$-coordinates.

\smallskip
\noindent
{\bf Step~5.} Now let us check that the choice of matrices in \eqref{key3, Mi} indeed satisfies the condition~\eqref{key2}. All we need to do is to evaluate the $(1,n+1)$-entry of 
\begin{align*}
\bbra I, \WW\kket_{\mw} =   \mw_{i_1}\left(e_{\ww_1}\right)\bullet \cdots \bullet \mw_{i_{n}}\left(e_{\ww_{n}}\right)\qquad \text{where } I \sim W.
\end{align*}
We refer to \eqref{notation, braket} for relevant notations, and to \eqref{two sums} and the ensuing remark in Step~2 for the reason of restricting to $I\sim W$.

This is completed by the two cases below.

\smallskip
\noindent
\underline{\emph{Case~1}: $\WW \neq I$.} There exists some $\ell \in \nset$ such that $i_\ell \neq \ww_\ell$. Thus $$\mw_{i_\ell}\left(e_{\ww_\ell}\right) = 0,$$ and hence $\bbra I, \WW\kket_{\mw}$ equals the zero matrix.

\smallskip
\noindent
\underline{\emph{Case~2}: $I=\WW$.} In this case we are interested in the term
\begin{align}\label{xx}
\bbra I, I\kket_{\mw} =   \mw_{i_1}\left(e_{i_1}\right)\bullet \cdots \bullet \mw_{i_{n}}\left(e_{i_{n}}\right),
\end{align}
which is a product of $n$ matrices.

For each $i \in \{i_1, \ldots, i_n\}$, by definition in  \eqref{key3, Mi} we have
\begin{equation}\label{yy}
\mw_i(e_i) := \sum_{t=1}^{r^W(i)} \left(E^{p^W(i|t)}_ {p^W(i|t)+1} - E^{p^W(i|t)+1}_{p^W(i|t)}\right).
\end{equation}
Thus
\begin{align}\label{linear comb}
\bbra I, I\kket_{\mw} &= \text{ linear combination of matrices of the form } E^{a_1}_{b_1} \bullet \cdots \bullet E^{a_n}_{b_n},
\end{align} 
where 
\begin{align*}
    \text{either $a_i + 1 = b_i$ or $a_i = b_i +1$ for each $i \in \nset$.}
\end{align*}

Now, thanks to the identity~\eqref{id on E} for elementary matrices, we have
\begin{align*}
E^{a_1}_{b_1} \bullet \cdots \bullet E^{a_n}_{b_n} = \delta_{b_1}^{a_2} \cdots \delta_{b_{n-1}}^{a_n}  E^{a_1}_{b_n}.
\end{align*}
If this term contributes nontrivially to the $(1,n+1)$-entry of $\bbra I,I\kket_{\mw}$, then one must have
\begin{align*}
a_1=1,\qquad b_n=n+1,\qquad \text{ and } a_{\ell+1}=b_\ell \text{ for each }\ell \in \{1,\ldots,n-1\}.
\end{align*}
This together with the triangle inequality and that  $|a_i-b_i|=1$ for each $i \in \nset$ (due to our choice of $a_i$, $b_i$) yields that
\begin{align*}
n &= b_n-a_1 \nonumber\\
&\leq |b_n-a_n| + |a_n-b_{n-1}| + |b_{n-1}-a_{n-1}| \nonumber\\
&\qquad \qquad + |a_{n-1}-b_{n-2}| + \ldots + |a_2-b_1|+|b_1-a_1|\nonumber\\
&= \sum_{i=1}^n |b_i-a_i| = n,
\end{align*}
so equalities must hold everywhere. But this is only possible when
\begin{align*}
(a_1, \ldots, a_n) = (1,\ldots, n-1) \qquad \text{and}\qquad (b_1, \ldots, b_n) = (2, \ldots, n). 
\end{align*}

Therefore, one infers from \eqref{linear comb} and the previous paragraph that the only nontrivial contribution to the $(1, n+1)$-entry of $\bbra I, I\kket_{\mw}$ comes from a single term of the form
\begin{align}\label{lambda}
\lambda E^1_2 \bullet E^2_3 \bullet \cdots\bullet E_n^{n-1} \qquad \text{for some $\lambda \in \R$.}
\end{align} 

Recall the expression \eqref{xx} for $\bbra I, I\kket_{\mw}$; each term on the right-hand side of which is expressed via \eqref{yy}. The term in \eqref{lambda} arises only if the following condition holds:
\begin{align*}
\text{For each $\ell \in \nset$, there is $t \in \left\{1, \ldots, r^W(i_\ell)\right\}$ such that $p^{W}\left(i_\ell|t\right)=\ell$}.
\end{align*}
In view of Definition~\ref{def: r and p}, this condition is tantamount to the following:
\begin{align*}
\text{For each $\ell \in \nset$, the index $i_\ell$ must occur at location $\ell$ in the word $W$}.
\end{align*}
This clearly implies $W=(i_1, \ldots, i_n)=I$.

Also, by an  inspection of \eqref{xx} and \eqref{yy}, the coefficient $\lambda$ in \eqref{lambda} must be 1. Thus we have verified \eqref{key2}, which completes the proof of Theorem~\ref{thm: main}.  \end{proof}

\begin{remark} \label{inv method remark} 

A careful examination reveals that the term $-E^{p^W(i|t)+1}_{p^W(i|t)}$ in our expression~\eqref{key3, Mi} for $\mw_i$ is somewhat immaterial to the proof  of Theorem~\ref{thm: main}. We include it only to ensure that $\mw_i$ takes value in the space $\so(k,\R)$ of skew-symmetric matrices. This leaves open the possibility for adapting our proof  to Cartan developments into more general Lie groups. 

The choice of $\mw_i$ is non-unique. Indeed, if we allow $\mw_i$ to take values in $\su(k)$ rather than restricting to $\so(k,\R)$, we may choose for instance $\2 E^{p^{W}(i|t)}_{p^{W}(i|t)} + E^{p^{W}(i|t)}_{p^{W}(i|t) + 1} - E^{p^{W}(i|t) + 1}_{p^{W}(i|t)}$ in lieu of $E^{p^{W}(i|t)}_{p^{W}(i|t) + 1} - E^{p^{W}(i|t) + 1}_{p^{W}(i|t)}$ in \eqref{key3, Mi}.
\end{remark}

\section{{Determination of measures over tensors of non-zero ROC}\label{sec: ESig}}

We now discuss implications of Theorem~\ref{thm: main} to the determination of Radon measures on subspaces of the extended tensor algebra $\tens$ and on the space of finite-$p$-variation paths.

\subsection{Determination of measures on \texorpdfstring{$\E_R$}{}} \label{Sub: determination_measures}

Recall from \eqref{t, def} the space of tridiagonal skew-symmetric matrices $\mathfrak{t}(k,\R)$. Also recall the definitions of the ROC of a tensor $X \in \tens$ and of $\E_R$ (see Definition~\ref{def: ROC}). The following version of the Stone--Weierstrass Theorem can be found, \emph{e.g.}, in \cite[Exercise~7.14.79]{measure}:

\begin{lemma}\label{lemma: SW}
Let $\mu$ and $\nu$ be Radon measures on a topological space $E$, and let $\mathcal{F}$ be a family of bounded continuous functions that separates points over $E$, contains the constant function $1$, and is closed under pointwise multiplication. If $\int_E f\,\dd\mu = \int_E f\,\dd\nu$ for every $f \in \mathcal{F}$, then $\mu = \nu$.
\end{lemma}

Fix any $R > 0$. As in \cite{cl}, we equip $\E_R$ with any Hausdorff topology such that each $\widetilde{M} \in \art$ is continuous. We have shown in Theorem~\ref{thm: main} that 
\begin{align*}
    \art :=  \left\{\widetilde{M} : \E_R \longrightarrow \bigoplus_{k \geq 1} \mathfrak{gl}(k;\R):\, M \in \bigoplus_{k \geq 1} \mathcal{L}\left(\R^d; \mathfrak{t}(k,\R)\right),\, ||M||_{\rm op} \leq R \right\}
\end{align*}
separates points over $\E_R$. That is, for any $X_1 \neq X_2$ in $\E_R$, there exists some $\widetilde{M_0} \in \art$ such that $\widetilde{M_0} (X_1) \neq \widetilde{M_0} (X_2)$. Note that $\mathfrak{t}(k;\R)$ is \emph{not} an algebra: this motivates us to consider
\begin{align*}
    \clart := \text{ the smallest algebra containing } \art.
\end{align*}
By an abuse of notation, we identify matrix algebras with the algebras of continuous functions of the matrix coefficients (denoted as $\mathcal{C}$ in \cite[p.4058, \S 4]{cl}). 

The result below follows immediately from Lemma~\ref{lemma: SW}. It together with Proposition~\ref{propn: ROC} (3)  directly implies \cite[Corollary~4.12]{cl}, which states that the Borel probability measures on the subspace of group-like elements $\mathcal{G}\left(\R^d\right)\subset\tens$ are determined by $\mathcal{A}$ defined in \eqref{A, def, new}.

\begin{proposition} \label{prop: E_R_measure}
Let $R > 0$ and $\mathcal{B} \subseteq \E_R$ be a Borel set such that every $f \in \clart$ is bounded on $\mathcal{B}$. Suppose $\mu$ and $\nu$ are Radon measures on $\E_R$. If $\int_\mathcal{B} f\,\dd\mu = \int_\mathcal{B} f\,\dd\nu$ for every $f \in \clart$, then $\mu = \nu$ on $\mathcal{B}$.
\end{proposition}

\begin{remark} \label{rmk: E_infty}
\hfill
\begin{enumerate}
    \item
    In Proposition~\ref{prop: E_R_measure}, any compact subset of $\E_R$ is a valid choice of $\mathcal{B}$.
    \item 
    For the case $R = \infty$, $\E_\infty$ can be topologised as a Polish space (see \cite[Corollary~2.4]{cl}). Hence, every finite Borel measure on $\E_\infty$ is Radon.
    \item 
    $\E_R$ is a Polish space under the metric $ d_R(X, Y) := \sum_{n=0}^\infty \|\pi_n(X) - \pi_n(Y)\| R^n$ for any $0<R<\infty$.
    \item 
    The space of group-like elements $\mathcal{G} \left(\R^d\right)$ is a closed subspace of $\E_\infty$. The arguments in \cite[\S4]{cl} lead to the characterisation: $\clart\Big|_{\mathcal{G} \left(\R^d\right)} = {\bf Clos}\left[\art\big|_{\mathcal{G} \left(\R^d\right)}\right]$. 
\end{enumerate}
\end{remark}

\subsection{Determination of measures on path spaces} \label{Sub: determine_path_measure}

Fix $1 \leq p < 2$. Denote by $\mathcal{V}_p\left([0, T], \R^d\right)$ the space of bounded-$p$-variation paths in $\R^d$ over $[0, T]$. Also denote its time-augmentation by
\[
\tilde{\mathcal{V}}_p\left([0, T], \R^d\right) := \left\{ (\gamma_t, t)_{t \in [0, T]} : \gamma \in \mathcal{V}_{p}\left([0, T], \R^d\right)\right\}.
\]
Clearly $\tilde{\mathcal{V}}_p\left([0, T], \R^d\right) \subset \mathcal{V}_p\left([0, T], \R^{d + 1}\right)$. Also, the signature maps $\tilde{\mathcal{V}}_p\left([0, T], \R^d\right)$ injectively into $\mathcal{G} \left(\R^d\right)$, which is a closed subspace of $\E_\infty$ (see Remark~\ref{rmk: E_infty}). 

Recall the definition of path development $\gamma \mapsto \mathscr{D}_\gamma$ from \eqref{dev, def}. As a direct consequence of Theorem~\ref{thm: main}, the path development separates points over $\tilde{\mathcal{V}}_p\left([0, T], \R^d\right)$.

\begin{corollary}\label{cor: new}
    For any $\gamma_1 \neq \gamma_2$ in $\tilde{\mathcal{V}}_p \left([0, T], \R^d\right)$, there exists $k \geq 1$ and $M_0 \in \lin \left(\R^d; \mathfrak{t}(k, \R)\right)$ such that $\mathscr{D}_{\gamma_1}(M_0) \neq \mathscr{D}_{\gamma_2}(M_0)$. Specifically, for any $n \geq 1$ such that $\pi_n \big[S(\gamma_1)\big] \neq \pi_n \big[S(\gamma_2)\big]$, we can find one such $M_0$ in $\mathcal{L}\left(\R^d; \mathfrak{t}(n + 1, \R)\right)$.
\end{corollary}

For ease of notation, write 
\begin{equation*}
    \mathscr{D}_\gamma^{[\mathfrak{g}]} := \text{the restriction of the path development $\mathscr{D}_\gamma$ to $\mathcal{L}\left(\R^d; \mathfrak{g}=\bigoplus_{k=1}^\infty \mathfrak{g}(k)\right)$.}
\end{equation*}
Here $\mathfrak{g}$ is a matrix Lie algebra graded by $k \in \mathbb{N}$ such that $\mathfrak{g}(k) \subset \mathfrak{gl}(k;\C)$. In this work, we are interested in the following three choices of $\mathfrak{g}(k)$: 
$$\mathfrak{t}(k, \R) \subset \mathfrak{so}(k) \subset \mathfrak{su}(k).$$

As in \S\ref{Sub: determination_measures} above, we equip $\tilde{\mathcal{V}}_p\left([0, T], \R^d\right)$ with any Hausdorff topology such that for any $M \in \bigoplus_{k \geq 1} \mathcal{L}\left(\R^d; \mathfrak{su}(k)\right)$, the unitary path development $\mathscr{D}_\cdot (M): \tilde{\mathcal{V}}_p\left([0, T], \R^d\right) \rightarrow \bigoplus_{k \geq 1} \mathcal{U}(k)$ is continuous. One valid choice is to equip $\tilde{\mathcal{V}}_p\left([0, T], \R^d\right)$ with the $p$-variation norm \cite{book}, which makes it an inseparable Banach space. The continuity of $\mathscr{D}_\cdot(M)$ holds because the development of a path is the unique solution of a linear differential equation driven by the path; see Remark~\ref{rmk:CDE}. The same argument holds if we replace $\mathfrak{su}(k)$ and $\mathcal{U}(k)$ by $\mathfrak{so}(k)$ and $\mathcal{O}(k)$, respectively.

Note that for $\mathfrak{g} \in \{\mathfrak{su}, \mathfrak{so}\}$, $\mathscr{D}^{[\mathfrak{g}]}_\cdot$ is a set of bounded continuous functions on $\tilde{\mathcal{V}}_p\left([0, T], \R^d\right)$ that contains the constants. In addition, for any $M_1 \in \mathcal{L}\left(\R^d; \mathfrak{g}(k_1)\right)$ and $M_2 \in \mathcal{L}\left(\R^d; \mathfrak{g}(k_2)\right)$, 
\begin{align} \label{eq: group_prod}
    \mathscr{D}^{[\mathfrak{g}]}_\cdot (M_1) \otimes \mathscr{D}^{[\mathfrak{g}]}_\cdot (M_2) = \mathscr{D}^{[\mathfrak{g}]}_\cdot (M_3)\qquad \text{ on $\tilde{\mathcal{V}}_p\left([0, T], \R^d\right)$},
\end{align}
where $M_3 = M_1 \otimes I_{k_2} + I_{k_1} \otimes M_2 \in \mathcal{L}\left(\R^d; \mathfrak{g}(k_1 k_2)\right)$. Therefore, from the Stone-Weierstrass Theorem (Lemma~\ref{lemma: SW}) we deduce that
\begin{proposition}\label{prop: determine_path_measure}
Continuous $\R$-valued functions on any compact subspace of $\tilde{\mathcal{V}}_p\left([0, T], \R^d\right)$ can be uniformly approximated by path developments $\mathscr{D}_\cdot ^{\mathfrak{g}}\,(\mathfrak{g} \in \{\mathfrak{su}, \mathfrak{so}\})$.

Meanwhile, let $\mu$ and $\nu$ be Radon measures on $\tilde{\mathcal{V}}_p\left([0, T], \R^d\right)$ with $1\leq p <2$. If every function in  $\mathscr{D}_\cdot ^{\mathfrak{g}}\,(\mathfrak{g} \in \{\mathfrak{su}, \mathfrak{so}\})$ has equal integrals with respect to $\mu$ and $\nu$, then $\mu = \nu$.
\end{proposition}

Proposition~\ref{prop: determine_path_measure} does not cover the case $\mathfrak{g} = \bigoplus_{k\geq 1}\mathfrak{t}(k;\R)$ since it is not an algebra.

Now we specialise to the expected signature of a path. If the ROC $\rho(\mathbb{E}(S(\gamma)))=\infty$, then 
\begin{equation*}
    \mathbb{E} \big[\mathscr{D}_\gamma (M)\big] = \mathbb{E} \left[\widetilde{M}(S(\gamma)) \right] = \widetilde{M} \big[\mathbb{E} (S(\gamma)) \big].
\end{equation*}
By Proposition~\ref{prop: determine_path_measure}, the expected signature uniquely determines Radon measures on $\tilde{\mathcal{V}}_p\left([0, T], \R^d\right)$ with $1\leq p <2$. Also, in view of \S~\ref{sec:generating_func}, $\widetilde{M} \big[\mathbb{E} (S(\gamma)) \big]$ is the generating function of the expected signature. Hence, the expected signature can be recovered from
\[
\left\{ 
\mathbb{E} \big[\mathscr{D}_\gamma (M)\big] : M \in \bigoplus_{k \geq 1} \mathcal{L}\left(\R^d; \mathfrak{t}(k, \R)\right) 
\right\}
\]
by the same algorithm in Theorem~\ref{thm: main}. We thus arrive at the following:
\begin{proposition} \label{prop: rpcf}
Let $\mu$ and $\nu$ be Radon measures on $\tilde{\mathcal{V}}_p\left([0, T], \R^d\right)$; $1\leq p<2$. Assume that their expected signatures exist and are of infinite ROC. The following are equivalent:
\begin{itemize}
    \item 
    $\mu = \nu$;
    \item 
    $\mathbb{E}_{\gamma \sim \mu} \big[\mathscr{D}_\gamma(M)\big] =  \mathbb{E}_{\gamma^\prime \sim \nu} \big[\mathscr{D}_{\gamma^\prime}(M)\big]$ for all $M \in \bigoplus_{k \geq 1} \lin\left(\R^{d}; \mathfrak{t}(k, \R)\right)$.
\end{itemize}
\end{proposition}

Under the condition that the signature transform 
\begin{align*}
    S: \tilde{\mathcal{V}}_p\left([0, T], \R^d\right) \longrightarrow \mathcal{S}_p := S\left(\tilde{\mathcal{V}}_p\left([0, T], \R^d\right)\right) \subset \mathcal{G} \left(\R^d\right)
\end{align*}
is a homeomorphism onto its image, one may easily obtain an analogous result to Proposition~\ref{prop: determine_path_measure} on the determination of Radon measures on $\mathcal{S}_p$.

\begin{remark} \label{rmk: geometric}
    In Proposition~\ref{prop: rpcf},  if one only considers measures on the time augmentation of the path space $\mathcal{V}_{0, p} \subset \mathcal{V}_p$, then every finite Borel measure is Radon as $\mathcal{V}_{0, p}$ is a Polish space  (\emph{cf}. \cite[Proposition~5.36]{friz2010multidimensional}).     $\mathcal{V}_{0, p}$ is the closure of smooth paths under the $p$-variation norm, and is known as the space of geometric $p$-rough paths. See Friz--Victoir \cite[Definition~9.15]{friz2010multidimensional}. 

\end{remark}

\section{Numerical experiments on hypothesis testing}\label{sec: ht}

Now we present various numerical results based on the theories developed earlier in this paper, especially \S\ref{sec: ESig}. In \S\ref{sec:pcfd}, we propose a new distance function for measures on path spaces. Then, in \S\ref{sec: ht_sub}, we apply this proposed distance function on the task of hypothesis testing (\emph{i.e.}, deciding whether two stochastic processes have the same distribution on the path space). The numerical results demonstrate that our proposed method significantly outperforms the signature-based MMD \cite{chevyrev2022signature}, and achieves a comparable performance to PCFD-based approaches \cite{ours}, while offering reduced time complexity by an order of magnitude. 

For readers interested in reproducing the numerical experiments, the implementation code is available in the GitHub repository:  \href{https://github.com/DeepIntoStreams/UnitaryDevInversion}{\color{blue}{https://github.com/DeepIntoStreams/UnitaryDevInversion}}.

\subsection{A new distance of stochastic processes}  \label{sec:pcfd}

Let $1 \leq p < 2$. In the sequel, we consider probability measures on $\tilde{\mathcal{V}}_{0, p}$, \emph{i.e.}, the time augmentation of the space of geometric $p$-rough paths, which is a Polish space (Remark~\ref{rmk: geometric}). Denote  
\begin{align*}
\mathcal{P}_\infty \left(\tilde{\mathcal{V}}_{0, p}\right) &:= \bigg\{ \text{Borel probability measures on $\tilde{\mathcal{V}}_{0, p}$ }\\
&\qquad \qquad \text{ whose expected signatures exist and with ${\rm ROC}=\infty$}\bigg\}.
\end{align*}
Then, by Proposition~\ref{prop: rpcf} and Remark~\ref{rmk: geometric} \emph{cf.} \cite{ours}), the Path Characteristic Function (PCF)
\[
\Phi_\mu (M) := \mathbb{E}_{\gamma \sim \mu} \big[\mathscr{D}_\gamma (M)\big] = \mathbb{E}_{\gamma \sim \mu} \left[\widetilde{M}(S(\gamma))\right]
\]
fully characterises $\mu \in \mathcal{P}_\infty \left(\tilde{\mathcal{V}}_{0, p}\right)$, as long as its domain of definition contains $ \bigoplus_{k \geq 1} \lin\left(\R^{d}; \mathfrak{t}(k, \R)\right)$.  This motivates us to introduce the following 

\begin{definition}\label{def_RPCFD} 
Let $\mathcal{M}:=(\mathcal{M}_k)_{k \geq 1}$, where $\mathcal{M}_k$ is a probability measure on $\mathcal{L}\left(\R^d; \mathfrak{t}(k, \R)\right)$. For any $\mu, \nu \in \mathcal{P}_\infty \left(\tilde{\mathcal{V}}_{0, p}\right)$, their Restricted Path Characteristic Function Distance (RPCFD) with respect to $\mathcal{M}$ is defined by
\begin{eqnarray*}
{\rm RPCFD}_{\mathcal{M}}^2(\mu, \nu):= \sup_{k \geq 1} \mathbb{E}_{M \sim \mathcal{M}_k} \Big[\big|\big| \Phi_\mu(M)- \Phi_\nu(M) \big|\big|_{\rm HS}^2\Big],
\end{eqnarray*}
where $||\cdot||_{\rm HS}$ denotes the matrix Hilbert--Schmidt norm.
\end{definition}

    In view of Proposition~\ref{prop: rpcf}, ${\rm RPCFD}_{\mathcal{M}}$ is a true distance as long as each $\mathcal{M}_k$ has a full support on $\mathcal{L}\left(\R^d; \mathfrak{t}(k, \R)\right)$. It implies that for any two different measures $\mu$ and $\nu$, there exists some $k \geq 1$ such that the $k$-th grade of ${\rm RPCFD}_{\mathcal{M}}^2(\mu, \nu)$  is non-zero; \emph{i.e.},
    \[
    \mathbb{E}_{M \sim \mathcal{M}_k} \Big[\big|\big| \Phi_\mu(M)- \Phi_\nu(M) \big|\big|_{\rm HS}^2\Big] \neq 0.
    \]
   Since $\mathfrak{t}(1, \R) \emb \mathfrak{t}(2, \R) \emb \cdots \emb \mathfrak{t}(k, \R) \emb \mathfrak{t}(k+1, \R)\cdots$, one may work with a sufficiently large grade $k$ and the corresponding expectation $\mathbb{E}_{M \sim \mathcal{M}_k}$ in numerical experiments. However, larger values of $k$ lead to an increased computational complexity and a higher risk of overfitting. Therefore, cross-validation is employed in practice to select the most appropriate value of $k$.

\begin{remark} \label{remark: OPCFD}
    As $\mathfrak{t}(k, \R) \subset \su(k)$, RPCFD can be viewed as a restricted version of PCFD \cite{ours} without losing the characteristicity of measures on the path space. Since 
    \begin{equation*}
        \dim(\mathfrak{t} (k, \R)) = k-1 < \dim(\su(k)) = k^2,
    \end{equation*} 
    RPCFD has a significant dimension reduction effect over PCFD. Likewise, based on Proposition~\ref{prop: determine_path_measure}, we define the Orthogonal Path Characteristic Function Distance (OPCFD) by restricting the domain of $M$ to $\mathcal{L}\left(\R^d; \bigoplus_{k=1}^\infty\so(k)\right)$. Notice that the characteristicity of OPCFD does not require the existence of expected signature. As $dim(\so(k)) = k(k-1)/2$, RPCFD reduces the dimension  more significantly than OPCFD, so we focus more on RPCFD. See Appendix~\ref{comparison} for numerical experiments on OPCFD.
\end{remark}

\begin{remark}
Typically, empirical time series data are mostly observed at discrete time points. As in \cite[Theorem~8.22 (i.3)]{friz2010multidimensional}, we approximate paths in $\tilde{\mathcal{V}}_{0, p}$ by interpolations (e.g., linear interpolation), and the corresponding Borel probability measures on $\tilde{\mathcal{V}}_{0, p}$ converge to the correct limit as the time step size tends to zero.  
\end{remark}

To enhance the discriminative power of RPCFD on given path distributions, given the grade $k \geq 1$, we propose the following trainable distance:
\begin{eqnarray*}
\sup_{\mathcal{M} \in \mathcal{P}(\mathcal{L}(\R^d; \mathfrak{t}(k, \R)))} \text{RPCFD}^2_{\mathcal{M}} (\mu, \nu).
\end{eqnarray*}
In practice, we approximate $\mathcal{M}$ by an empirical measure:
\begin{eqnarray*}
    \mathcal{M} \approx \frac{1}{K}\sum_{i = 1}^{K} \delta_{M_i},
\end{eqnarray*}
where $M_i \in \mathcal{L}\left(\R^d; \mathfrak{t}(k, \R)\right)$, $i \in \{1, \cdots, K\}$, for some $K \geq 1$. Note that the linear maps $\{M_i\}_{1 \leq i \leq K}$ can be fully parameterised. Therefore, the task of learning the optimal distribution $\mathcal{M}$ in \eqref{eqn_tr_dist} is transformed to learning the optimal parameters in the linear maps $\{M_i\}_{1 \leq i \leq K}$ which attain the supremum of the ``empirical RPCFD'':
\begin{eqnarray} \label{eqn_tr_dist}
    \sup_{\{M_i\}_{1 \leq i \leq K}} \frac{1}{K}\sum_{i = 1}^{K} \text{RPCFD}^2_{\delta_{M_i}} (\mu, \nu) = \sup_{\{M_i\}_{1 \leq i \leq K}} \frac{1}{K}\sum_{i = 1}^{K} || \Phi_{\mu}(M_{i}) - \Phi_{\nu}(M_{i})||_{\rm HS}^2.
\end{eqnarray}

\subsection{Numerical experiments on hypothesis testing} \label{sec: ht_sub}

We apply the trainable empirical RPCFD in \eqref{eqn_tr_dist} to the task of hypothesis testing on distributions of stochastic processes. More specifically, we first train the empirical RPCFD \eqref{eqn_tr_dist} on the training set. We then use the trained empirical RPCFD as the statistic to conduct a two-sample permutation test on the test set. Finally, we use the test power (the higher the better) and Type-I error (the lower the better) as the performance metrics.

\subsubsection*{\fbox{Two-sample permutation test}}

Suppose $X$ and $Y$ are two stochastic processes with distributions $\mathbb{P}_X$ and $\mathbb{P}_Y$, respectively. The hypotheses of the two-sample permutation test are:
\begin{eqnarray}
    && H_{0}: \mathbb{P}_{X} = \mathbb{P}_{Y} \quad \text{(null hypothesis)}, \nonumber \\
    && H_{1}: \mathbb{P}_{X} \neq \mathbb{P}_{Y} \quad \text{(alternative hypothesis)}. \nonumber
\end{eqnarray}
The test power is defined as the probability of correctly rejecting $H_0$ when $H_1$ is true, and the Type-I error denotes the probability of falsely rejecting $H_0$ when $H_0$ is true.

Denote by $\mathbf{X} = \{X_1, X_2, \cdots, X_m\}$ and $\mathbf{Y} = \{Y_1, Y_2, \cdots, Y_n\}$ two sets of i.i.d. samples drawn respectively from $\mathbb{P}_X$ and $\mathbb{P}_Y$. Given the statistic $T(\mathbf{X}, \mathbf{Y})$ (here $T$ is the trained empirical RPCFD), we apply all possible permutations $\sigma \in \Sigma_{m+n}$ to $\mathbf{Z}:=(\mathbf{X}, \mathbf{Y})$ and obtain
\begin{equation*}
    \mathcal{T} := \bigg\{T\big(\mathbf{Z}_{\sigma(1):\sigma(m)}, \mathbf{Z}_{\sigma(m+1):\sigma(m+n)}\big) : \; \sigma \in \Sigma_{m+ n}\bigg\}.
\end{equation*}
Under $H_0$, one can uniformly sample from $\mathcal{T}$ to approximate the distribution of $T(\mathbb{P}_X, \mathbb{P}_Y)$. Given the significance level $\alpha$, the null hypothesis $H_0$ will be rejected if $T(\mathbf{X}, \mathbf{Y})$ exceeds the $(1-\alpha)\%$ quantile of $\mathcal{T}$, and otherwise accepted.

The algorithm of our experiments is in Appendix~\ref{details}. See \emph{e.g.}, Chevyrev--Oberhauser \cite[\S7]{chevyrev2022signature} for more details on two-sample permutation tests.

\subsubsection*{Constructing training and test sets}
Following \cite[Example~B.12]{ours}, we set $X := \tilde{B}_{[0, T]}$ to be the time-augmented 3-dimensional Brownian Motion (BM), and $Y := \tilde{B}^{h}_{[0, T]}$ to be the time-augmented 3-dimensional fractional Brownian Motion (fBM) with the Hurst parameter $h \in ]0, 1[$. Note that $\mathbb{P}_{X} \neq \mathbb{P}_{Y}$ when $h \neq 0.5$. The training data are i.i.d. generated under $H_1$, and the test data (\emph{i.e.}, out-of-sample data) are i.i.d. generated under either $H_0$ or $H_1$, depending on the choice of performance metrics. See Appendix~\ref{details} for settings of hyperparameters.

\begin{remark}
For any $\e>0$, BM is of finite $(2 + \epsilon)$-variation and fBM is of finite $(1/h + \epsilon)$-variation. Neither is of finite 2-variation when $h \leq 0.5$.  Our numerics indicate that the trained empirical EPCFD is a powerful statistic to separate different measures on the highly oscillatory paths. 
\end{remark}

\subsubsection*{\fbox{RPCFD vs. PCFD}}

\begin{figure}[H]
    \centering
    \begin{minipage}{0.49\textwidth}
        \includegraphics[width=\textwidth]{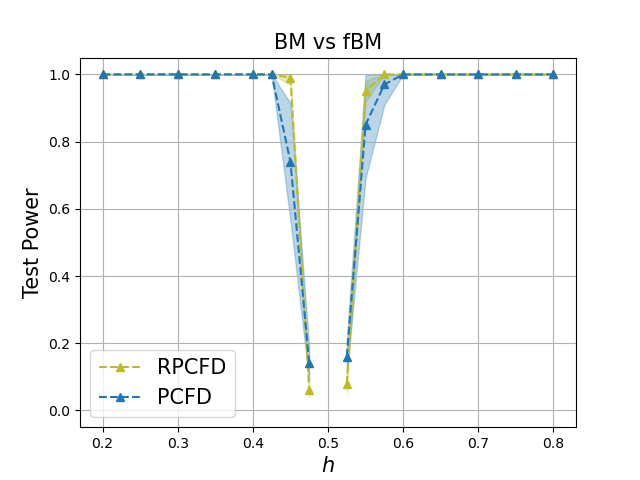}
    \end{minipage}
    \begin{minipage}{0.49\textwidth}
        \includegraphics[width=\textwidth]{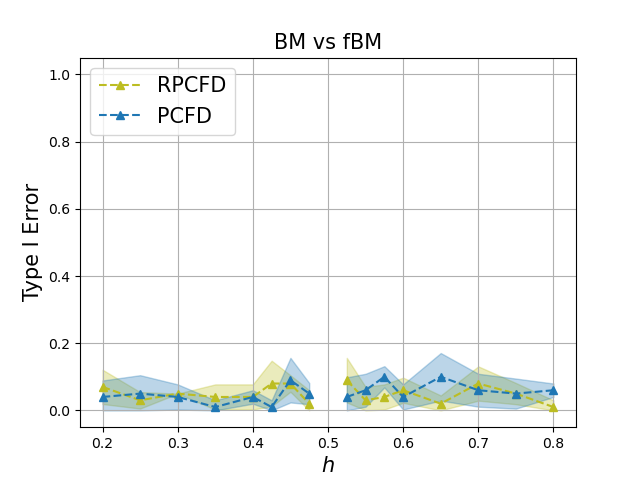}
    \end{minipage}
    
    \begin{minipage}{0.49\textwidth}
        \includegraphics[width=\textwidth]{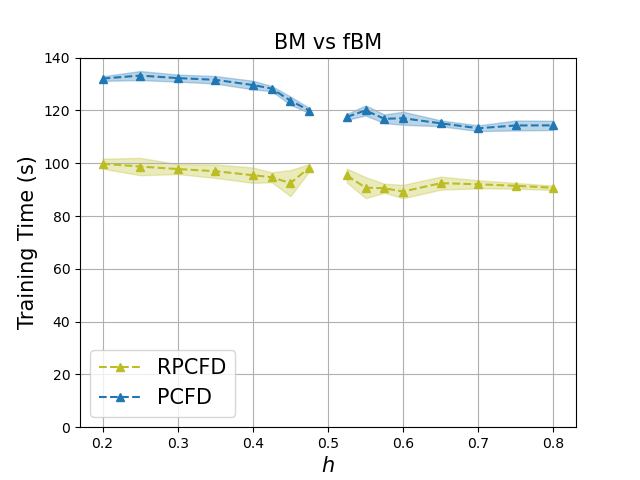}
    \end{minipage}
    \vspace{-3mm}
    \caption{(Upper Left) Test powers of the trained RPCFD and PCFD against various $h$, where the shaded area represents mean $\pm$ std over 5 runs. (Upper Right) Type-I errors of the trained RPCFD and PCFD. (Bottom) Training time of RPCFD and PCFD over 500 iterations. We fix $K = 8$ and $k = 5$ in all experiments.}
    \label{fig: RPCFD}
\end{figure}

We first compare the performance of our proposed trainable empirical RPCFD \eqref{eqn_tr_dist} with PCFD proposed in \cite{ours}. In our experiments, we set $K = 8$ and $k = 5$, and both distances are trained on the same training set by exactly the same method. As shown in Figure~\ref{fig: RPCFD} (Upper Left), the overall test power curves against various $h$ of RPCFD and PCFD are close. For $h = 0.45$ and $h = 0.55$, RPCFD achieves an obviously superior test power over PCFD. Figure~\ref{fig: RPCFD} (Upper Right) indicates that the two distances have comparable performances in Type-I error, while the bottom panel shows that the training time of RPCFD has been notably reduced ($\approx 22.9\%$ on average)  compared with that of PCFD. 

\begin{remark}
    The improved performance of RPCFD over PCFD indicates that restricting the linear maps to $\mathfrak{t}(k, \R)$ narrows down the search space and hence makes the learning more efficient. On the other hand, the training time of RPCFD is also reduced, since $\mathfrak{t}(k, \R)$ has only real entries and thus simplifies the computation (in contrast, all parameters in PCFD are complex) by a factor of $2$. In Figure~\ref{fig: OPCFD}, we repeat the experiments with RPCFD replaced by OPCFD (\emph{cf.} Remark~\ref{remark: OPCFD}) and report the experiment results. Similar results are observed, which supports our analysis.
\end{remark}

\begin{figure}[H]
    \centering
    \begin{minipage}{0.49\textwidth}
        \includegraphics[width=\textwidth]{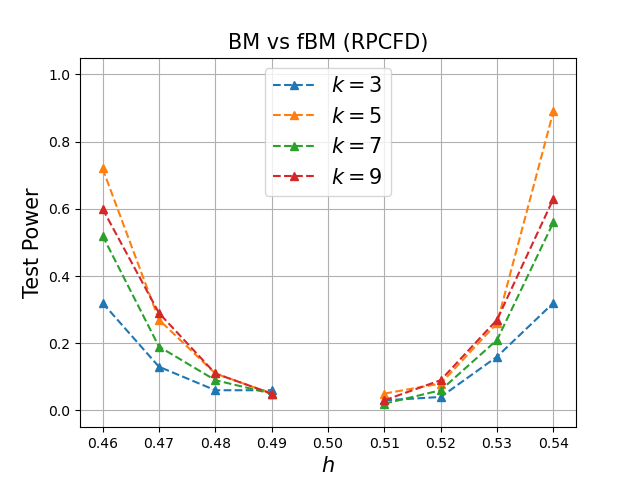}
    \end{minipage}
    \begin{minipage}{0.49\textwidth}
        \includegraphics[width=\textwidth]{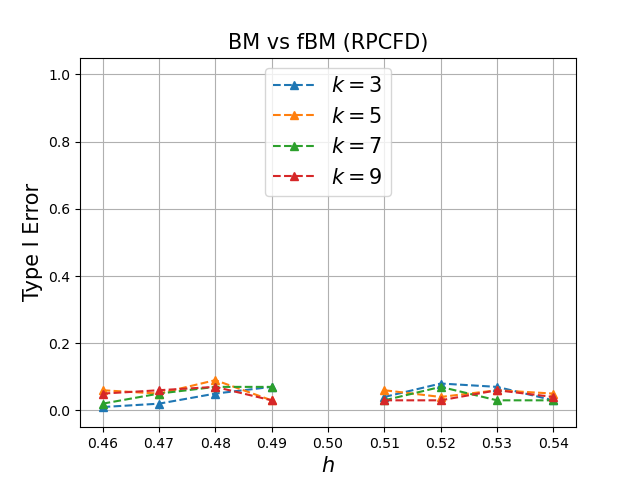}
    \end{minipage}
    \vspace{-3mm}
    \caption{Sensitivity analysis of RPCFD (Left: mean test power over 5 runs; Right: mean Type-I error over 5 runs). The $y$-axis and $x$-axis represent the power of the trained RPCFD and the Hurst parameter $h$, respectively. We fix $K = 8$ in all experiments.}
    \label{sensitivity_plot}
\end{figure}

The sensitivity analysis of RPCFD when $|h-0.5| \in [0.01, 0.04]$ is plotted in Figure~\ref{sensitivity_plot}. In general, increasing the matrix order $k$ improves the discriminative power. However, a larger value of $k$ may have the potential issues of overfitting and higher computational costs. Therefore, in practice, one should choose $k$ by hyperparameter tuning on a cross-validation set.

\subsubsection*{\fbox{Comparison with Sig-MMDs}}

In Appendix~\ref{comparison}, we provide a comprehensive summary of the test powers, Type-I errors, and time complexities of all the distances/models (\emph{i.e.}, RPCFD, OPCFD, PCFD and signature-based MMDs) in Tables~\ref{Tab:summary_table_type2}, \ref{Tab:summary_table_type1}, and \ref{Tab:summary_time}. According to Table~\ref{Tab:summary_table_type2}, the development-based distances (\emph{i.e.}, RPCFD, OPCFD, and PCFD) outperform the signature-based MMDs consistently for various $h$ in test power. In particular, when $h \in \{0.45, 0.55\}$, the development-based distances improve the test power from $[0.15, 0.17]$ to $[0.95, 0.99]$. Table~\ref{Tab:summary_table_type1} shows that the distances/models have close performances in Type-I error. In Table~\ref{Tab:summary_time}, we observe that the inference time of signature-based MMDs scales quadratically with respect to the sample size, while the dependency is linear for development-based methods (\emph{i.e.}, RPCFD, OPCFD, and PCFD). For example, when sample size $m=n=200$, the inference time of signature-based MMDs is over 10 times longer than that of the distance-based methods.

\medskip
\noindent
{\bf Acknowledgement}.
SL and HN are deeply indebted to Prof.~Haojie Ren for insightful discussions on hypothesis testing. All the authors thank Dr.~Lei Jiang for his help with running the codes on servers.

\medskip
\noindent
{\bf Funding sources}.
The research of SL is supported by NSFC Projects 12201399, 12331008, and 12411530065, Young Elite Scientists Sponsorship Program by CAST  2023QNRC001, the National Key Research $\&$ Development Programs 2023YFA1010900 and 2024YFA1014900, Shanghai Rising-Star
Program 24QA2703600, and the Shanghai Frontier Research Institute for Modern Analysis. HN is supported by the EPSRC under the program grant EP/S026347/1 and the Alan Turing Institute under the EPSRC grant EP/N510129/1. Both SL and HN are partially funded by the SJTU-UCL joint seed fund WH610160507/067 and the Royal Society International Exchanges 2023 grant (IEC/NSFC/233077).

\medskip
\noindent
{\bf Competing interests statement}. 
All the authors declare that there is no conflict of interest.

\newpage
\appendix
\section{Supplementary details of numerical experiments} \label{details}

\subsection{Algorithm of two-sample permutation test in \S\ref{sec: ht_sub}}
\hfill

\begin{algorithm}[H]
    \begin{algorithmic}[1]
        \Require $\alpha \in (0, 1)$: significance level; $N > 0$: \# of experiments; $M > 0$: \# of permutations; $\mathbf{X}$: samples from distribution $\mu$; $\mathbf{Y}$: samples from distribution $\nu$; $m > 0$: sample size of $\mathbf{X}$; $n > 0$: sample size of $\mathbf{Y}$; $T$: test statistic function; $H_0 \in \{1, 0\}$: whether the null hypothesis is true or false ($H_0 = 1$ if $\mu = \nu$, otherwise $H_0 = 0$).
        \State $\mathbf{Z} \gets$ Concatenate$(\mathbf{X}, \mathbf{Y})$
        \State $\text{num\_rejections} \gets 0$
        \State $i \gets 1$
        \While{$i \leq N$}
            \State $\mathcal{T} \gets$ EmptyList
            \State $j \gets 1$
            \While{$j \leq M$}
                \State $\sigma \sim $ Permutation($\{1, 2, \cdots, m + n\}$)
                \State $T_{\sigma} \gets T(\{\mathbf{Z}_{\sigma(1)}, \mathbf{Z}_{\sigma(2)}, \cdots, \mathbf{Z}_{\sigma(m)}\}, \{\mathbf{Z}_{\sigma(m + 1)}, \cdots, \mathbf{Z}_{\sigma(m + n)}\})$
                \State $\mathcal{T}\text{.append}(T_{\sigma})$
                \State $j \gets j + 1$
            \EndWhile
            \If{$T(\mathbf{X}, \mathbf{Y}) > (1-\alpha)\%$ quantile of $\mathcal{T}$}
                \State $\text{num\_rejections}\gets \text{num\_rejections} + 1$
            \EndIf
            \State $i \gets i + 1$
        \EndWhile
        \State $\text{ratio} \gets \text{num\_rejections }  / \text{ } N$
        \If{$H_0$}
            \State $\text{Type\_I\_error}\gets \text{ratio}$
            \State \Return $\text{Type\_I\_error}$
        \Else
            \State $\text{test\_power}\gets \text{ratio}$ 
            \State \Return $\text{test\_power}$
        \EndIf
    \end{algorithmic}
    \caption{Estimate the test power/Type-I error of the permutation test}
    \label{alg:power}
\end{algorithm}

\subsection{Hyperparameters and training procedure in \S\ref{sec: ht_sub}}

All experiments are implemented on a Tesla V100 GPU. The random seed of PyTorch and NumPy is set to be 0. When estimating the test power, 20 experiments are implemented (\emph{i.e.}, $N = 20$). In each experiment, 200 BM and 200 fBM trajectories are sampled as the out-of-sample test dataset (\emph{i.e.}, the construction of $\mathbf{X}$ and $\mathbf{Y}$ with $m = n = 200$), which is permuted 500 times to estimate the distribution of the permuted distance (\emph{i.e.}, $M = 500$). The significance level of the permutation test is set to be 5\% (\emph{i.e.}, $\alpha = 0.05$). 

When training each development-based distance (\emph{i.e.}, RPCFD, OPCFD and PCFD), 10000 BM and 10000 fBM trajectories (50 time steps on $[0, 1]$) are sampled as the training set. The training utilizes PyTorch's Adam optimiser with the $\beta$'s being $(0, 0.9)$. The fine-tuned learning rate of OPCFD and PCFD is 0.05 whereas that of RPCFD is 0.5. In each of the 500 iterations of the training, a mini-batch of 1024 trajectories is randomly selected from the training set to conduct the stochastic gradient ascent to maximise the corresponding empirical distance. After training, a pool of 10000 BM and 10000 fBM trajectories is independently sampled and used to conduct the permutation test. 

For the signature MMDs, we consider the linear and RBF kernels (see \cite[\S7.1]{chevyrev2022signature} for more details). For the RBF kernel, the fine-tuned $2 \sigma^2$ is set to be 0.1.

\section{Supplementary numerical results} \label{comparison}

\subsection*{OPCFD vs. PCFD}
\begin{figure}[H]
    \centering
    \begin{minipage}[b]{0.49\textwidth} 
        \includegraphics[width=\textwidth]{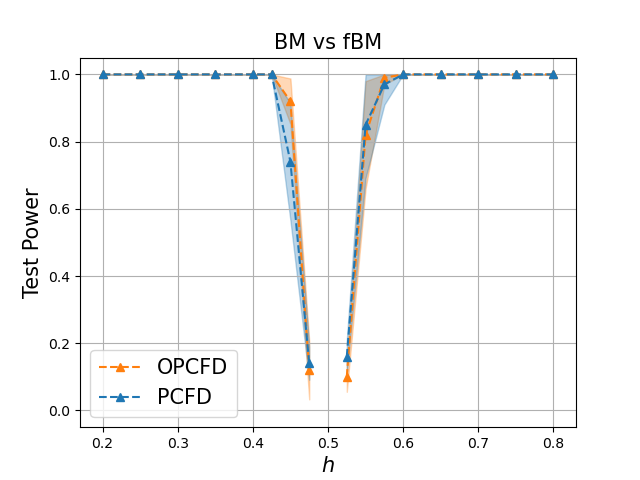}
    \end{minipage}
    \begin{minipage}[b]{0.49\textwidth} 
        \includegraphics[width=\textwidth]{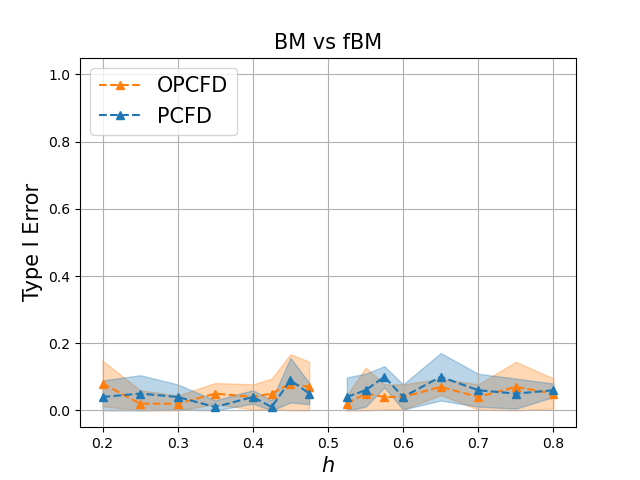}
    \end{minipage}
    \vspace{-2mm}
    \begin{minipage}[b]{0.49\textwidth} 
        \includegraphics[width=\textwidth]{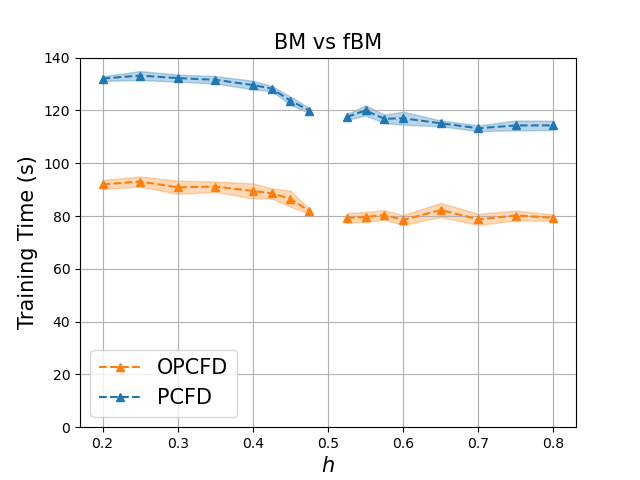}
    \end{minipage}
    \caption{(Upper Left) Test powers of the trained OPCFD and PCFD against various $h$, where the shaded area represents mean $\pm$ std over 5 runs. (Upper Right) Type-I errors of the trained OPCFD and PCFD. (Bottom) Training time of OPCFD and PCFD over 500 iterations (time reduction of OPCFD $\approx 31.0\%$ on average). We fix $K = 8$ and $k = 5$ in all experiments.}
    \label{fig: OPCFD}
\end{figure}

As shown in Figure~\ref{fig: OPCFD}, OPCFD shows a comparable performance to PCFD while achieving an average time reduction of $\approx 31.0\%$.

\newpage
\subsection*{Summarising tables}

Comprehensive tables for summarising the test power, Type-I error and time complexity of the distances/models involved in this paper are presented below.

\begin{table}[b]
\renewcommand{\arraystretch}{1.2}
\resizebox{0.75\columnwidth}{!}
{
\begin{tabular}{lccccccc}
\hline
& & \multicolumn{3}{c}{Developments} & & \multicolumn{2}{c}{Signature MMDs} \\ \cmidrule{3-5} \cmidrule{7-8}
$h$ & & RPCFD & OPCFD & PCFD & & Linear & RBF \\
\hline
0.2  & & $\textbf{1}\pm 0$ & $\textbf{1}\pm 0$ & $\textbf{1}\pm 0$ & & $\textbf{1}\pm 0$ & $\textbf{1}\pm 0$ \\
0.25 & & $\textbf{1}\pm 0$ & $\textbf{1}\pm 0$ & $\textbf{1}\pm 0$ & & $\textbf{1}\pm 0$ & $\textbf{1}\pm 0$ \\
0.3  & & $\textbf{1}\pm 0$ & $\textbf{1}\pm 0$ & $\textbf{1}\pm 0$ & & $0.9\pm 0.04$     & $\textbf{1}\pm 0$ \\
0.35 & & $\textbf{1}\pm 0$ & $\textbf{1}\pm 0$ & $\textbf{1}\pm 0$ & & $0.35\pm 0.07$    & $\textbf{1}\pm 0$ \\
0.4  & & $\textbf{1}\pm 0$ & $\textbf{1}\pm 0$ & $\textbf{1}\pm 0$ & & $0.09\pm 0.06$    & $0.97\pm0.03$    \\
0.425& & $\textbf{1}\pm 0$ & $\textbf{1}\pm 0$ & $\textbf{1}\pm 0$ & & $0.1\pm 0.05$     & $0.69\pm0.11$     \\
0.45 & & $\textbf{0.99}\pm 0.02$ & $0.92\pm 0.07$ & $0.74\pm 0.17$& & $0.04\pm 0.04$ & $0.15\pm0.05$   \\
0.475& & $0.06\pm 0.02$ & $0.12\pm 0.09$ & $\textbf{0.14}\pm 0.05$& & $0.01\pm 0.02$    & $0.04\pm0.02$  \\
0.525& & $0.08\pm 0.02$ & $0.1\pm 0.04$ & $\textbf{0.16}\pm 0.05$& & $0.05\pm0.02$    & $0.07\pm0.04$   \\
0.55 & & $\textbf{0.95}\pm 0.03$ & $0.82\pm 0.16$ & $0.85\pm 0.16$       & & $0.13\pm0.05$  &$0.17\pm0.04$   \\
0.575& & $\textbf{1}\pm 0$ & $0.99\pm 0.02$       & $0.97\pm 0.06$        & & $0.07\pm0.02$    &$0.5\pm0.10$  \\
0.6  & & $\textbf{1}\pm 0$ & $\textbf{1}\pm 0$ & $\textbf{1}\pm 0$ & & $0.05\pm0.03$ &$0.75\pm0.05$ \\
0.65 & & $\textbf{1}\pm 0$ & $\textbf{1}\pm 0$ & $\textbf{1}\pm 0$ & & $0.14\pm0.05$ & $\textbf{1}\pm 0$ \\
0.7  & & $\textbf{1}\pm 0$ & $\textbf{1}\pm 0$ & $\textbf{1}\pm 0$ & & $0.23\pm0.05$ & $\textbf{1}\pm 0$\\
0.75 & & $\textbf{1}\pm 0$ & $\textbf{1}\pm 0$ & $\textbf{1}\pm 0$ & &$0.28\pm0.07$ & $\textbf{1}\pm 0$ \\
0.8  & & $\textbf{1}\pm 0$ & $\textbf{1}\pm 0$ & $\textbf{1}\pm 0$ & &$0.46\pm0.07$ & $\textbf{1}\pm 0$\\
\hline
\end{tabular}}
\vspace{3mm}
\caption{Test power of the distances when $h \neq 0.5$ in the form of mean $\pm$ std over 5 runs. For RPCFD, OPCFD and PCFD, fix $K = 8$ and $k = 5$. For the RBF signature MMD, fix $2 \sigma^2 = 0.1$.}
\label{Tab:summary_table_type2}
\renewcommand{\arraystretch}{1}
\end{table}

\begin{table}[b]
\renewcommand{\arraystretch}{1.2}
\resizebox{0.75\columnwidth}{!}{
\begin{tabular}{lccccccc}
\hline
& & \multicolumn{3}{c}{Developments} & & \multicolumn{2}{c}{Signature MMDs} \\ \cmidrule{3-5} \cmidrule{7-8}
$h$ & & RPCFD & OPCFD & PCFD & & Linear & RBF  \\
\hline
0.2  &&$0.07\pm0.05$ & $0.08\pm0.07$ & $0.04\pm0.05$ && $0.04\pm0.09$  &  $0.04\pm0.04$ \\
0.25 &&$0.03\pm0.02$ & $0.02\pm0.04$ & $0.05\pm0.05$ && $0.03\pm0.03$ &  $0.04\pm0.04$ \\
0.3  &&$0.05\pm0$    & $0.02\pm0.02$ & $0.04\pm0.04$ && $0.04\pm0.04$  &  $0.01\pm0.02$ \\
0.35 &&$0.04\pm0.04$ & $0.05\pm0.03$ & $0.01\pm0.02$ && $0.03\pm0.04$  &  $0.07\pm0.07$ \\
0.4  &&$0.04\pm0.04$ & $0.04\pm0.04$ & $0.04\pm0.02$ && $0.06\pm0.05$ &  $0.04\pm0.04$ \\
0.425&&$0.08\pm0.07$ & $0.05\pm0.04$ & $0.01\pm0.02$ && $0.03\pm0.03$  &  $0.04\pm0.04$ \\
0.45 &&$0.08\pm0.02$ & $0.08\pm0.09$  & $0.09\pm0.07$ && $0.05\pm0.04$ & $0.02\pm0.03$  \\
0.475&&$0.02\pm0.04$ & $0.07\pm0.07$ & $0.05\pm0.03$ && $0.05\pm0.04$  &  $0.07\pm0.06$ \\
0.525&&$0.09\pm0.07$ & $0.02\pm0.02$ & $0.04\pm0.06$ && $0.03\pm0.03$  &  $0.04\pm0.02$ \\
0.55 &&$0.03\pm0.04$ & $0.05\pm0.08$ & $0.06\pm0.05$ && $0.07\pm0.06$  & $0.05\pm0.05$  \\
0.575&&$0.04\pm0.04$ & $0.04\pm0.04$ & $0.1\pm0.03$ && $0.02\pm0.03$  &  $0.06\pm0.07$ \\
0.6  &&$0.06\pm0.04$ & $0.04\pm0.04$ & $0.04\pm0.04$ && $0.05\pm0.04$  & $0.05\pm0.04$   \\
0.65 &&$0.02\pm0.02$ & $0.07\pm0.02$ & $0.1\pm0.07$ && $0.07\pm0.06$  & $0.03\pm0.04$   \\
0.7  &&$0.08\pm0.05$ & $0.04\pm0.04$ & $0.06\pm0.05$ && $0.05\pm0.05$  & $0.03\pm0.04$   \\
0.75 &&$0.05\pm0.03$ & $0.07\pm0.07$ & $0.05\pm0.04$ && $0.06\pm0.07$  & $0.04\pm0.07$   \\
0.8  &&$0.01\pm0.02$ & $0.05\pm0.04$ & $0.06\pm0.02$ && $0.07\pm0.07$  & $0.04\pm0.04$   \\
\hline
\end{tabular}}
\vspace{3mm}
\caption{Type-I error of the distances when $h \neq 0.5$ in the form of mean $\pm$ std over 5 runs. For RPCFD, OPCFD and PCFD, fix $K = 8$ and $k = 5$. For the RBF signature MMD, fix $2 \sigma^2 = 0.1$.}
\label{Tab:summary_table_type1}
\renewcommand{\arraystretch}{1}
\end{table}

\begin{table}[t]
\renewcommand{\arraystretch}{1.2}
\resizebox{\columnwidth}{!}{
\begin{tabular}{cccccccc}
\hline
& & \multicolumn{3}{c}{Developments} & & \multicolumn{2}{c}{Signature MMDs} \\ \cmidrule{3-5} \cmidrule{7-8}
 & & RPCFD & OPCFD & PCFD & & Linear & RBF  \\
\hline
\multicolumn{1}{c|}{$m=n$}& \multicolumn{7}{c}{Inference time (seconds)}\\
\midrule
\multicolumn{1}{c|}{$10$}  &&$175.35\pm8.24$ & $168.26\pm10.34$ & $167.69\pm2.20$ && $95.12\pm 0.21$  &  $122.9\pm 0.32$ \\
\multicolumn{1}{c|}{$50$} &&$212.54\pm8.78$ & $204.02\pm3.68$ & $236.21\pm5.46$ && $402.69\pm0.23$ &  $533.43\pm 0.33$ \\
\multicolumn{1}{c|}{$100$} &&$333.08\pm9.23$ & $319.57\pm3.01$ & $376.69\pm4.25$ && $1329.73\pm0.57$ &  $1760.18\pm0.29$ \\
\multicolumn{1}{c|}{$200$} &&$398.87\pm22.3$ & $400.63\pm26.06$ & $510.71\pm20.08$ && $5257.11\pm0.98$ &  $6958.81\pm1.33$ \\
\midrule
\multicolumn{1}{c|}{Mini-batch size}& \multicolumn{7}{c}{Training time (seconds over 500 iterations)}\\
\midrule 
\multicolumn{1}{c|}{1024}  &&$92.46\pm4.93$ & $86.53\pm3.07$ & $123.71\pm1.56$ && $-$  &  $-$ \\
\hline
\end{tabular}}
\vspace{3mm}
\caption{Inference time of the permutation test across different sample sizes and the training time of RPCFD, OPCFD and PCFD before conducting the permutation test. The result is in the form of mean $\pm$ std over 5 runs. For RPCFD, OPCFD and PCFD, fix $K = 8$ and $k = 5$. For the RBF signature MMD, fix $2 \sigma^2 = 0.1$. We fix the Hurst parameter of fBM to be 0.45 in all experiments.}
\label{Tab:summary_time}
\renewcommand{\arraystretch}{1}
\end{table}


\begin{thebibliography}{99}
\bibitem{bdmn}
H. Boedihardjo, J. Diehl, M. Mezzarobba, H. Ni, The expected signature of Brownian motion stopped on the boundary of a circle has finite radius of convergence, \textit{Bull. Lond. Math. Soc.} \textbf{53}(1) (2021) 285--299. https://doi.org/10.1112/blms.12420

\bibitem{bg}
H. Boedihardjo, X. Geng, ${\rm SL}_2(\R)$-developments and signature asymptotics for planar paths with bounded variation, \textit{ArXiv preprint} (2020) ArXiv: 2009.13082; accepted by \textit{Rev. Mat. Iberoam.} \textbf{39} (2023) 1973--2006. https://doi.org/10.4171/RMI/1397

\bibitem{c1}
K.-T. Chen, Iterated integrals and exponential homomorphisms, \textit{Proc. Lond. Math. Soc.} \textbf{4}(3) (1954) 502--512. https://doi.org/10.1112/plms/s3-4.1.502

\bibitem{c2}
K.-T. Chen, Integration of paths, geometric invariants and a generalized Baker--Hausdorff formula, \textit{Ann. of Math.} \textbf{65}(1) (1957) 163--178. https://doi.org/10.2307/1969671

\bibitem{c3}
K.-T. Chen, Integration of paths -- a faithful representation of paths by noncommutative formal power series, \textit{Trans. Amer. Math. Soc.} \textbf{89}(2) (1958) 395--407. https://doi.org/10.1090/S0002-9947-1958-0106258-0

\bibitem{ck}
I. Chevyrev, A. Kormilitzin, A primer on the signature method in machine learning, \textit{ArXiv preprint} (2016), ArXiv:1603.03788. 
https://doi.org/10.48550/arXiv.1603.03788

\bibitem{cl}
I. Chevyrev, T. Lyons, Characteristic functions of measures on geometric rough paths, \textit{Ann. Probab.} \textbf{44}(6) (2016) 4049--4082. https://doi.org/10.1214/15-AOP1068


\bibitem{alg1}
A. Giambruno, A. Valenti, On minimal $\ast$-identities of matrices, \textit{Linear $\&$ Multilinear Algebra} \textbf{39} (1995) 309--323. https://doi.org/10.1080/03081089508818405

\bibitem{alg2}
A. Giambruno, M. Zaicev, Polynomial identities and asymptotic methods, 
\textit{Math. Surveys Monogr.} \textbf{122} (2005). https://doi.org/10.1090/surv/122

\bibitem{hl}
B.~M. Hambly, T.~J. Lyons, Uniqueness for the signature of a path of bounded variation and the reduced path group, \textit{Ann. of Math.} \textbf{171}(1) (2010) 109--167. https://doi.org/10.4007/annals.2010.171.109

\bibitem{unique}
H. Boedihardjo, X. Geng, T.~J. Lyons, D. Yang, The signature of a rough path: uniqueness, \textit{Adv. Math.} \textbf{293} (2016) 720--737. https://doi.org/10.1016/j.aim.2016.02.011

\bibitem{heyer}
H. Heyer, Probability Measures on Locally Compact Groups, Springer, Berlin, 1977.

\bibitem{lln}
D. Levin, T.~J. Lyons, H. Ni, Learning from the past, predicting the statistics for the future, learning an evolving system, \textit{ArXiv preprint} (2013), ArXiv: 1309.0260. https://doi.org/10.48550/arXiv.1309.0260

\bibitem{lini}
S. Li, H. Ni, Expected signature of stopped Brownian motions on $d$-dimensional $C^{2,\alpha}$-domains has finite radius of convergence everywhere: $2 \leq d \leq 8$, \textit{J. Funct. Anal.} \textbf{282}(12) (2022). https://doi.org/10.1016/j.jfa.2022.109447

\bibitem{ours'}
H. Lou, S. Li, H. Ni, Path development network with finite-dimensional Lie group representation, \textit{ArXiv preprint} (2022), ArXiv: 2204.00740. 
https://doi.org/10.48550/arXiv.2204.00740

\bibitem{ours}
H. Lou, S. Li, H. Ni, PCF-GAN: generating sequential data via the characteristic function of measures on the path space, \textit{ArXiv preprint} (2023), ArXiv: 2305.12511; accepted by \emph{NeurIPS 2023}.

\bibitem{icm}
T.~J. Lyons, Rough paths, signatures and the modelling of functions on streams, \emph{Proceedings of the International Congress of Mathematicians}, Korea, 2014.

\bibitem{book}
T.~J. Lyons, M. Caruana, T. L\'{e}vy, Differential Equations Driven by Rough Paths, Springer-Verlag Berlin, Heidelberg, 2007.

\bibitem{friz2010multidimensional}
P.~K. Friz, N.~B. Victoir, Multidimensional stochastic processes as rough paths: theory and applications, Cambridge University Press, 2010.

\bibitem{phd}
T.~J. Lyons, H. Ni, Expected signature of Brownian motion up to the first exit time from a bounded domain, \textit{Ann. Probab.} \textbf{43}(5) (2015) 2729--2762. https://doi.org/10.1214/14-AOP949

\bibitem{lx}
T.~J. Lyons, W. Xu, Hyperbolic development and inversion of signature, \textit{J. Funct. Anal.} \textbf{272}(7) (2017) 2933--2955. https://doi.org/10.1016/j.jfa.2016.12.024

\bibitem{nixu}
H. Ni, W. Xu, Concentration and exact convergence rates for expected Brownian signatures, \textit{Electron. Commun. Probab.} \textbf{20} (2015) 1--11. https://doi.org/10.1214/ECP.v20-3636

\bibitem{hurst}
R. Passeggeri, On the signature and cubature of the fractional Brownian motion for $h>1/2$, \textit{Stoch. Process. Appl.} \textbf{130}(3) (2020) 1226--1257. https://doi.org/10.1016/j.spa.2019.04.013

\bibitem{esig}
T. Fawcett, Problems in stochastic analysis. Connections between rough paths and non-commutative harmonic analysis, \textit{PhD diss.} University of Oxford (2003).

\bibitem{morrill2019signature}
J. Morrill, A. Kormilitzin, A. Nevado-Holgado, S. Swaminathan, S. Howison, T.~J. Lyons, The signature-based model for early detection of sepsis from electronic health records in the intensive care unit, \textit{2019 Computing in Cardiology}, Singapore, 2019, pp. 1--4. https://doi.org/10.22489/CinC.2019.014

\bibitem{ni2023conditional}
H. Ni, L. Szpruch, M. Wiese, S. Liao, B. Xiao, Conditional Sig-Wasserstein GANs for Time Series Generation, \textit{Math. Finance} (Special Issue on Machine Learning in Finance) \textbf{34}(2) (2024), 622--670. https://doi.org/10.1111/mafi.12423

\bibitem{kidger2019deep}
P. Bonnier, P. Kidger, I.~P. Arribas, C. Salvi, T.~J. Lyons, Deep signature transforms, \textit{NeurIPS 2019}, Article No.279, 3105--3115.

\bibitem{chevyrev2022signature}
I. Chevyrev, H. Oberhauser, Signature moments to characterize laws of stochastic processes, \textit{J. Mach. Learning Res.} \textbf{23}(1) (2022) 7928--7969. 

\bibitem{measure}
V.I. Bogachev, Measure theory, Springer-Verlag Berlin, Heidelberg, 2007.

\end{thebibliography}
\end{document}